\documentclass[11pt]{article}

\usepackage{amsmath,amssymb,amsthm, latexsym,amsfonts,epsfig, color, geometry, enumitem, hyperref,float,multirow,bm,array,xfrac,graphicx}
\usepackage{xcolor,colortbl}
\usepackage[linesnumbered,ruled,vlined]{algorithm2e}
\newcommand{\ip}[2]{\langle #1 , #2 \rangle}    

\newcommand{\st}{\operatorname{s.t.}}

\newcommand{\cB}{{\mathcal B}}

\newcommand{\A}{{\mathcal A}}

\newcommand{\cC}{{\mathcal C}}

\newcommand{\Sym}{\mathbb S}

\newcommand\ran{\mathcal R}
\newcommand\R{{\mathbb R}}

\newcommand{\tra}{\text{tr\hspace{0.05cm}}}
\newcommand{\tr}{^{\top}}

\newcommand{\eps}{\varepsilon}

\newtheorem{theorem}{Theorem}

\newtheorem{proposition}{Proposition}
\newtheorem{corollary}{Corollary}
\newtheorem{definition}{Definition}
\newtheorem{remark}{Remark}

\newtheorem{assumption}{Assumption}
\newtheorem{problem}{Problem}

\definecolor{mygreen}{cmyk}{0.82,0.11,1,0.25}

\newcommand{\new}[1]{\textcolor{black}{ #1}}

\usepackage{endnotes}
\let\footnote=\endnote

\title{Two-Stage Robust Quadratic Optimization with Equalities and its Application to Optimal Power Flow}

 \author{
 Olga Kuryatnikova\thanks{Econometric Institute, Erasmus University Rotterdam, the Netherlands, \href{mailto:kuryatnikova@ese.eur.nl}{kuryatnikova@ese.eur.nl}, the corresponding author} \and
Bissan Ghaddar\thanks{Ivey Business School, Western University, \href{mailto:bghaddar@ivey.ca}{bghaddar@ivey.ca}}
\and Daniel K. Molzahn\thanks{School of Electrical and Computer Engineering, Georgia Institute of Technology, \href{mailto:molzahn@gatech.edu}{molzahn@gatech.edu}}
}
\usepackage{amsopn}

\definecolor{Gray}{gray}{0.95}
\newcolumntype{g}{>{\columncolor{Gray}}c}

\begin{document}
\maketitle
\begin{abstract}
In this work, we consider two-stage quadratic optimization problems under ellipsoidal uncertainty. In the first stage, one needs to decide upon the values of a subset of optimization variables (control variables). In the second stage, the uncertainty is revealed and the rest of the optimization variables (state variables) are set up as a solution to a known system of possibly non-linear equations. This type of problem occurs, for instance, in optimization for dynamical systems, such as electric power systems as well as gas and water networks.  We propose a  {convergent} iterative algorithm to build a sequence of approximately robustly feasible solutions with an improving objective value. At each iteration, the algorithm optimizes over a subset of the feasible set and uses affine approximations of the second-stage equations while preserving the non-linearity of other constraints. We implement our approach and demonstrate its performance on M{\sc atpower} instances of AC Optimal Power Flow. This paper focuses on quadratic problems, but the approach is suitable for more general setups.
\end{abstract}


%


\section{Introduction} \label{sec:intro}

In many optimization problems, data is not completely known in advance. One of the main approaches to deal with this lack of information is robust optimization (RO). It  assumes that the data lies in a predefined set of scenarios and the constraints have to be satisfied for any realization of the data in that set. RO does not require any knowledge about the distribution of the uncertain data and is useful when the distribution is hard to estimate and feasibility for a certain set of parameters is important. In reality, some optimization variables represent decisions that must be made before the actual realization of the uncertain data while other variables can be adjusted after the uncertain data becomes known. To account for such situations in  RO,  \emph{two-stage adjustable} robust optimization (ARO) was introduced in \cite{BenTal2004}. Two-stage ARO problem includes two types of variables: the first-stage variables that are fixed and the second-stage variables that may change depending on the uncertainty realization. We refer to these variables as control and state variables, respectively, since this setting is typical for optimal control problems. ARO gives rise to more flexible decisions than robust optimization and thus could be less conservative.  

 {Several approaches have been developed to solve ARO problems,  most of them are approximations since solving ARO to optimality can be NP-hard even if the original problem without uncertainty is a linear program  \cite{BenTal2004}. The hardness comes from defining the relations between the first- and second-stage decision variables, called \emph{decision rules}, which are usually not specified and must be optimized.  In certain problems occurring in practice, the exact functional form of the second-stage rules is predetermined but unknown since it is defined \emph{implicitly} via, for instance, a system of non-linear equalities that are challenging, or even impossible, to eliminate. Implicitly defined decision rules are rarely considered in the literature while they occur in many applications, some of which are described below. This paper aims to partly close this gap and studies general quadratic ARO with decision rules implicitly defined by a system of quadratic equalities. }

With rapidly increasing uncertainties in both the demand and supply of resources, networked infrastructure problems are important applications of ARO with implicitly defined decision rules. Examples of relevant networked infrastructures include electricity (see~\cite{ChanceBienstock}), natural gas (see~\cite{JuanGasNetw, misra2020}), and water (see~\cite{stuhlmacher2020}). Operators must ensure that these systems remain in acceptable states despite uncertainties while also considering performance criteria such as operating costs. ARO provides a way to balance these potentially competing concerns, as discussed in~\cite{AdjROSurvey} and \cite{misra2020}.  As an illustrative application of ARO, we have chosen the so-called \emph{robust AC optimal power flow} (ACOPF) problem that provides minimum cost operating points for electric power systems. ACOPF can be formulated as a non-convex quadratic optimization problem, as discussed in~\cite{GhaddarOPF}. Even in the absence of uncertainties, ACOPF is NP-Hard (see~\cite{bienstock2015nphard}), and solving such problems under uncertainty for instances of realistic size is still a challenge,  {which we explain  in detail in Section~\ref{subsec:acopf}.}


\subsection{Contribution} \label{subsec:contr}  

 {This paper presents two main contributions. First, we propose a convergent iterative solution approach for two-stage non-convex quadratic problems under uncertainty. The resulting solutions are feasible for the underlying problem without uncertainty (in other words, the \emph{nominal problem}) and approximately  feasible for the ARO version of the problem. Second, we implement the proposed approach for ACOPF, and it outperforms the benchmark methods.} \medskip \\ 
\noindent     \textbf{Solution framework for quadratic ARO under ellipsoidal uncertainty with implicit decision rules defined by quadratic equality constraints.} The equalities increase the difficulty of ARO. Therefore, we address equality and inequality constraints separately. We construct piecewise affine approximations of the implicit decision rules. In particular, we express the state variables as functions of the control and uncertainty variables using the first-order Taylor approximations of the original implicit rules on small subsets. To our knowledge, we present the first algorithm for ARO which uses implicit second-stage rules.\\
 For each piece of the piecewise affine approximation, we eliminate the equalities and second-stage variables and obtain a standard non-linear quadratic problem in the  \emph{first-stage variables} under ellipsoidal uncertainty. This problem is reformulated into a  {semidefinite program (SDP)} with (possibly) quadratic constraints.   We suggest using an alternating projections algorithm  to find locally optimal solutions for this SDP in the presence of quadratic constraints.  {Such a first-stage solution is  feasible for the nominal problem, and possible constraint violations under uncertainty are limited by the parameters of our method.} \medskip \\ 
\noindent     \textbf{Implementation for ACOPF and comparison with benchmarks.}  We apply the proposed framework to ACOPF with uncertainty in power supply and demand. The ACOPF is formulated as a non-convex quadratically constrained quadratic problem, and an ellipsoidal uncertainty set is considered. As shown in Section~\ref{subsec:acopf}, our numerical results for ACOPF demonstrate the effectiveness of our approach in comparison to other approaches from the literature on small to moderate-size instances ranging from  6 to 118 buses, with the potential to consider larger instances. 


\subsection{Existing ARO solution approaches} \label{subsec:literature}
 Several approaches have been developed to solve ARO problems,  { most of them are approximations since solving ARO to optimality can be NP-hard even if the original problem without uncertainty is an LP  \cite{BenTal2004}.}
 The hardness comes from defining the relations between the first- and second-stage decision variables, called \emph{decision rules}, which are usually not specified and must be optimized. 
One of the most popular ways to solve general ARO is an approximation where the second-stage decision variables can be written as affine functions of the uncertain parameters. This approach is 
proposed by \cite{BenTal2004}, and its current state-of-the art is discussed in \cite{georghiou2020primal}.   {Besides affine decision rules, one could consider piecewise constant decision rules. Such rules are constructed by partitioning the uncertainty set into subsets and implementing a fixed second-stage decision for each subset. Naturally,  the two types of decision rules could be combined into piecewise affine decision rules, which is suggested in \cite{splittingUnc}. Fixing the form of decision rules could restrict flexibility, so one could miss the optimal ARO solution. As an alternative, there exist convergent relaxations which gradually add constraint violating uncertain scenarios to the problem \cite{Bienstock2008,Bertsimas2012,RowAndColumnGen}. These approaches work best if the uncertainty set is polyhedral and the nominal problem can be solved efficiently, which is different from the setup in our paper.   \\
Our approach has features of both described ARO frameworks. On the one hand, we use piecewise affine decision rules. However, those cannot be arbitrary good rules, as they must reflect an actual process represented via a system of equalities to prevent second-stage infeasibilities. Moreover, the equalities are not solvable for some first-stage decisions and uncertainty values. Thus, our method constructs piecewise rules that are close to the real ones when those exist and indicates a lack of second-stage solutions otherwise. Using approximations might still lead to second-stage infeasibilities, but one can control those by increasing the number of affine pieces. This behavior resembles the approaches ensuring feasibility on a subset of the uncertainty realizations.}

Finally, let us discuss the literature addressing non-linear ARO with implicitly defined decision rules. 
Recently, some progress has been made in methods suitable for such ARO, and we outline the differences between our approach and these methods. 
First, we do not consider any assumptions on convexity and concavity, which makes our work different from all results in convex and linear ARO.  {Another advantage is that our approach, while using SDP, is not based on SDP approximations of the original problem, such as in~\cite{JuanGasNetw,lasserreQuant,AdaptRobACOPF,venzke2017,Xie18}. We allow keeping some non-linearities from the original problem as opposed to classical SDP approaches that linearize all constraints. Additionally, independently of how much non-linearity one keeps, the resulting SDP constraints are of the size of the number of uncertain parameters.} On the contrary, other approaches use SDPs of at least the size of the number of first- and second-stage variables. The number of uncertain parameters is often substantially smaller than the number of variables. For example, in ACOPF the number of variables depends on the number of buses in the system,  uncertainty usually occurs in loads or renewable generators, and not all buses have those. Hence, the size of SDP constraints in our approach could be substantially smaller than in the above-mentioned papers.

Our approach differs from other robust optimization  methods where non-linear constraints are linearized, such as~\cite{RobustACOPF}. We linearize equalities only, linearize locally, and our approximations are closely related to the original constraints via Taylor series. The results in this paper are also distinct from~\cite{ACOPFRobFeas} since in the latter work a robust solution is obtained by iteratively tightening the inequality constraints. In~\cite{GenNLRobust}, the authors  tackle general non-linear optimization problems but use an alternative ARO formulation and thus a distinct solution strategy. Our approach is close in spirit to the approach in \cite{roald2018} which uses local Taylor expansions too. 
However, the authors of the latter paper work with chance constraints, use full linearizations and address specifically ACOPF. Similarly, the approach in~\cite{biefel2022affinely} is tailored to linear complementarity problems.

This paper is organized as follows. In Section~\ref{sec:problem}, we present the problem formulation and motivate our solution approach. In Section~\ref{sec:Eq}, we describe the proposed algorithms in detail. 
 In Section~\ref{subsec:acopf}, we evaluate our approach on ACOPF instances. Finally, in Section~\ref{sec:aviol} we present conclusions  and discuss directions for future work.

\section{Problem formulation and general framework} \label{sec:problem}
We begin this section with the notation. We denote the  range of matrix $A$ by $\ran(A)$. We denote the space  of  $n{\times}n$  symmetric matrices  by $\Sym^n$ and for $A,B \in \Sym^n$, the trace inner product of $A$ and $B$ is denoted by  $\ip{A}{B}:={\rm trace}(AB)$. 
We use the notation $[n]$ for the set $\{1,\dots,n\}$. For a vector $V$ of length $n$, we denote the $i^{th}$ entry of $V$ by $V_i$. We say that two continuous maps $f$ and $g$ on a compact set $A\subset \R^n$  are \emph{$\eps$-close to each other} if $\sup_{x\in A} \|f(x)-g(x)\| \le \eps$ for some given norm $\|\cdot\|$. We say that two vectors (or elements of a vector space) $a,b$ are $\eps$-close to each other if  $ \|a-b\| \le \eps$ for some properly defined norm.  

In this paper, we deal with \emph{quadratically constrained quadratic problems (QCQP)}. Our approach generalizes to problems of higher degree as we explain later, but we focus on QCQP for simplicity. Now, let $n_\zeta$, $n_y$, $n_x$, $m_{eq}$, and $m_{in}$ be natural numbers and consider semialgebraic sets $S_y\in \R^{n_y}$, $S_x\in \R^{n_x}$ defined by quadratic constraints. Consider \emph{quadratic} mappings $f: \R^{n_y} \to \R$, $G:\R^{n_\zeta} \times \R^{n_x} \times \R^{n_y} \to \R^{m_{in}}$, $L_x: \R^{n_\zeta} \times \R^{n_x} \to \R^{m_{eq}}$ and $L_y: \R^{n_\zeta} \times \R^{n_y} \to \R^{m_{eq}}$. Now, we can introduce the \emph{two-stage nominal} problem.
\begin{problem}[Nominal two-stage quadratic optimization problem]
\begingroup
\linespread{0.6}\selectfont{} \label{pr:nom} 
\begin{align*}
\inf_{y,x} \quad &   f(y) \nonumber  \\
 \quad \st \quad &  y \in S_{y}\nonumber   \\ 
&  L_i(y,x)=0 \text{ for all } i\in [m_{eq}]\\
&  G_i(y,x)\ge 0 \text{ for all } i\in [m_{in}]  \\
& x \in S_x. 
\end{align*}
\addtocounter{equation}{-1} 
\endgroup
\end{problem}  
Problem~\eqref{pr:nom} occurs frequently in dynamic systems optimization, which inspired this paper. In such problems, an operator sets up the values of control variables $y$, and the state variables $x$ are determined afterwards according to a system of equations defining the equilibrium. Some examples of dynamic systems optimization are energy problems (ACOPF as described in Section~\ref{subsec:acopf}, optimal power dispatch presented by \cite{bingane2019tight}), water problems (the valve placement problem  by \cite{GhaddarWater}) or gas problems (passive gas network feasibility problem studied by \cite{JuanGasNetw}). Moreover, problem~\eqref{pr:nom} describes the more general class of \emph{bilevel optimization problems} where after setting the values of $y$ in the first stage, the second stage variables are chosen from the set of optimal solutions of the second-stage optimization problem.  The KKT optimality conditions of the second-stage problem can be written as a system of equalities, where the final set of second-stage variables $x$ consists of the original second-stage variables and Lagrange multipliers of the second-stage problem. As a result, one obtains the so-called \emph{complementarity formulation} of the initial bilevel optimization problem, which has the form of problem~\eqref{pr:nom} if the initial problem was linear. A typical example of a bilevel optimization problem is the Stackelberg competition in economics. Some examples of bilevel  problems in engineering can be found in~\cite{raghunathan2003math,baumrucker2010mpec}. More information about linear complementarity problems and robust optimization approaches for them can be found in~\cite{biefel2022affinely}. 
\begin{remark} \label{rem:bin}
Problem~\eqref{pr:nom} allows for binary variables and absolute values. One can write binary constraints on a variable $a$ as quadratic equality constraints $a=a^2$. One can write the constraint $a=|b|$ as $a^2=b^2, a\ge 0$. 
\end{remark}
An example of problem~\eqref{pr:nom} that is the main use case in this paper is the ACOPF problem. It can be written as follows using the general notation above.
\begingroup
\begin{subequations}
\linespread{0.6}\selectfont{}
\begin{align*}
 \min_{y,x} \quad &   y\tr Py+ p\tr y +p_0 \ \nonumber  \\
 \quad \st \quad &  A y \le b  \nonumber   \\   
&  x\tr Q_i x  + r_i=y_i \ \text{ for all } i\in [m_{eq}] \ \\
 &  x\tr Q_j x + r_j \ge 0 \ \ \text{ for all } j\in [m_{in}],  
\end{align*}
\end{subequations}
\addtocounter{equation}{-1}
\endgroup
where $y$ are active powers and voltage magnitudes on PV buses, and $x$ are voltages in rectangular form, see Section~\ref{subsec:acopf} for the full problem formulation. Another example is the valve setting problem in water distribution networks from \cite{GhaddarWater} with the following formulation:
\begingroup
\begin{subequations}
\linespread{0.6}\selectfont{}
\begin{align*}
 \min_{y,x} \quad &   d\tr y \ \nonumber  \\
 \quad \st \quad &  A_1 y \le b_1, A_2 x  \le b_2  \\ 
 &  x\tr Q_i x + q_i\tr x + c_i\tr y+ r_i =0 \ \text{ for all } i\in [m_{eq}] \ \\
 &  x\tr Q_j x + q_j\tr x + y\tr C_j y + c_j\tr y + x\tr R_j y+ r_j \le 0 \ \ \text{ for all } j\in [m_{in}],  
\end{align*}
\end{subequations}
\addtocounter{equation}{-1}
\endgroup
where $y$ are pressure heads and valve placement indicators, and $x$ are flow rates and absolute values of flow rates. 

In practice, given the uncertain demand and supply that is encountered in such applications, one often aims at solving problem~\eqref{pr:nom}  under uncertainty, which results in the next formulation. 
\begingroup
\def\theequation{\theproblem}
\begin{problem}[Two-stage quadratic ARO problem]  \begin{subequations}\label{pr:0gen_po} 
\linespread{0.6}\selectfont{}
\begin{align}
z= \inf_{y} \quad &   f(y) \nonumber  \\
 \quad \st \quad &  y \in S_{y}\nonumber   \\ 
& \text{and for any } \zeta \in \Omega \text{ there exists }  x \text{ such that the following holds:} \nonumber \\
&  L_i(y,\zeta,x)=0 \text{ for all } i\in [m_{eq}]\label{eq:L0}\\
&  G_i(y,\zeta,x)\ge 0 \text{ for all } i\in [m_{in}] \label{ineq:G0} \\
& x \in S_x. \label{ineq:x0}
\end{align}
\end{subequations} 
\end{problem} 
\endgroup

Problem~\eqref{pr:0gen_po} is a two-stage ARO problem with the \emph{uncertain parameter}~$\zeta$. The first stage happens before the uncertainty realization. At this stage, one assigns values to the \emph{control} variables $y \in S_y$. The second stage happens after the uncertainty realization. At this stage, one has to choose the best feasible value of the \emph{state} variables $x$ for the given uncertainty realization. The goal is to select a value for $y \in S_y$ such that there would be a feasible solution $x$ in the second stage for any uncertainty realization  $\zeta \in \Omega$. 
Any solution to problem~\eqref{pr:0gen_po} is feasible for the underlying nominal problem and robust against potential uncertainty. We obtain the nominal problem by setting $\zeta=0$ in problem~\eqref{pr:0gen_po}. Next, we summarize main assumptions used in this paper.  
\begin{assumption}[The characteristics of the objective and the constraints] \label{as:basic}
\begin{enumerate}[label=(\alph*)]
\item $S_x$, $S_y$ are compact sets. \label{as:shapeSet} 
\item All functions in the equality constraints~\eqref{eq:L0} are continuously differentiable. \label{as:shapeZeta0}
\item  All functions besides the equality constraints~\eqref{eq:L0} are polynomials of degree at most two. \label{as:shapeZeta} \
\emph{This assumption is for simplicity. In fact, the problem only has to be quadratic in $\zeta$, the requirements for other variables are milder. Our approach would be still applicable if the problem were polynomial of higher degree in $y$; see Section~\ref{sec:Ineq} for more details. The inequality constraints could be general polynomials in $x$ as well. If the inequalities have higher degree in $x$, we can obtain degree-two polynomials using variable substitution and increasing the number of state variables and equality constraints. We emphasize that this procedure does not influence the final size of the problems we solve since this size only depends on the number of control and uncertainty variables, i.e., $y$ and~$\zeta$.}
\item $\Omega$ is an ellipsoidal uncertainty set of the form:
\begin{align}
\Omega=&\left \{\zeta\in \R^{n_\zeta}: \zeta \tr \Sigma \zeta + \sigma \tr \zeta + r  \ge 0, \ j\in [m_{\zeta}]\right \}, \label{def:OmegaQuad} 
\end{align}
where $\Sigma$ is negative semidefinite.  {$\Omega$ has a non-empty interior.}  \label{as:shapeOmega} \\
\emph{We need assumptions on the degree of $\zeta$ and the shape of $\Omega$ to efficiently eliminate $\zeta$ from the problem using the S-lemma~\cite{Yakubovich}; see Section~\ref{subsec:uncert}. We have chosen an ellipsoidal uncertainty set as the base case since it is convenient for our approach; frequently appears in the literature as being less conservative than, for instance, box uncertainty; and has interpretations from both robust and chance-constrained perspectives; see, e.g.,~\cite{golestaneh2018ellipsoidal,chen2010cvar}.}
\end{enumerate}
\end{assumption}

\begin{assumption}[Assumptions without loss of generality] \label{as:wlog}
\begin{enumerate}[label=(\alph*)]
\item $S_x$ is defined by inequalities.  \label{as:sx} \\ 
 \emph{The assumption is w.l.o.g. since if the definition of  $S_x$ contains equalities, they could be moved to~\eqref{eq:L0} as a preprocessing step. If $S_x$ becomes unbounded afterwards, a large  ball constraint for $x$ can be added to the problem.}
  \item   $n_x = m_{eq}$ and there are no redundant equality constraints in~\eqref{eq:L0}. \label{as:basic4}   \\
\emph{The assumption is w.l.o.g.  If $n_x\le m_{eq}$ and there are redundant constraints, they could be detected and eliminated as a preprocessing step. If $n_x > m_{eq}$, then some state variables are free and could be added to the pool of control variables.}
\item Objective $f$ is amenable for optimization (e.g., convex quadratic or linear).  \label{as:fconv}\\
\emph{The assumption is w.l.o.g. since if $f$ is not convex, we introduce the epigraph control variable and add the epigraph constraint to $S_y$.}
\end{enumerate}
\end{assumption}

Our goal is to approximate the original problem by \emph{a problem in control variables} $y$ only, since they represent the actual decisions to implement. The following idea motivates our approach: if we could analytically solve equalities~\eqref{eq:L0} for the second-stage variables $x$, we would express $x$ as a function of $y$ and $\zeta$ obtaining the second-stage decision rule $p: \R^{n_y \times n_\zeta}\to \R^{n_x}$. Then problem~\eqref{pr:0gen_po} would be equivalent to the following problem:
\begingroup
\begin{problem} \label{pr:0gen_poSubs} 
\begin{subequations} 
\linespread{0.6}\selectfont{}
\begin{align}
z= \inf_{y} \quad &   f(y) \nonumber  \\
 \quad \st \hspace{0.3cm} &  y \in S_{y}\nonumber   \\
 &  G_i(y,\zeta,p(y,\zeta))\ge 0 \text{ for all } \zeta \in \Omega, \ i\in [m_{in}] \nonumber \\
&  p(y,\zeta) \in S_x \text{ for all } \zeta \in \Omega. \nonumber
\end{align}
\end{subequations} 
\end{problem}\endgroup
That is, if a known second-stage decision rule exists, then we could substitute it in the problem and eliminate the equalities and state variables to obtain a classical (not adjustable) robust optimization problem.  Clearly, if the second-stage variables are determined from a system of non-linear equalities, there might be no unique analytical expression for the decision rule in the problem. However, such expressions exist on small subsets of $S_y \times S_\zeta$ under known conditions according to the implicit function theorem (see, e.g., \cite{Spivak}). In essence, the approach we suggest restricts the optimization problem to such  subsets and replaces the implicit decision rule by its first-order Taylor approximation on each subset. We present the approach in detail in the next sections.

\section{Piecewise affine approximations of equality constraints} \label{sec:Eq}

Our approach looks at piecewise affine approximations of the second-stage decision rules $x=p(y,\zeta)$ in problem~\eqref{pr:0gen_poSubs}, which are given implicitly via the system of equalities~\eqref{eq:L0}, using Taylor series. In addition to its simplicity, the main advantages of the Taylor approximation are its good fit for the original function around the approximation point  {and the possibility to construct it for an implicit function. }
For the construction, we need some Jacobians of the system of equalities~\eqref{eq:L0}. Denote $ L(y,\zeta,x):=\left[\begin{array}{cc}
       L_1 (y,\zeta,x) \\
       \dots \\
        L_{m_{eq}}(y,\zeta,x) \end{array} \right]$.
        For an $(\hat{y}, \hat{\zeta}, \hat{x}) \in S_y\times \Omega \times S_x$,  let
$J^{\hat y, \hat \zeta, \hat x}_{y}$,  $J^{\hat y, \hat \zeta, \hat x}_{\zeta}$ and $J^{\hat y, \hat \zeta, \hat x}_{x} $ be the Jacobians of $L$ with respect to $y,\zeta$, and $x$, respectively, evaluated at  $(\hat y, \hat \zeta, \hat x)$. By Assumption~\ref{as:basic}\ref{as:basic4}, the Jacobian $J^{\hat y, \hat \zeta, \hat x}_{x} $ is a square matrix of the size $n_{eq} \times n_{eq}$. 
 \begin{theorem}[The implicit function theorem: \cite{Spivak} Theorem 2-12] \label{thm:ift} Let $L(y,\zeta, x):  {\R^{n_y}} \times   {\R^{n_\zeta}}  \times \R^{n_x} \to \R^{n_{eq}}$ be a continuously differentiable function. Let $(\hat{y}, \hat{\zeta}, \hat{x}) \in {\R^{n_y}} \times   {\R^{n_\zeta}}  \times \R^{n_x}$ be such that $L(\hat{y}, \hat{\zeta}, \hat{x})=0$ and the Jacobian $J^{\hat y, \hat \zeta, \hat x}_{x} $  of $L$ with respect to $x$ at $(\hat{y}, \hat{\zeta}, \hat{x})$ is non-singular. Then there exists an open set $U \subseteq  {\R^{n_y}} \times   {\R^{n_\zeta}} $ with $(\hat y, \hat \zeta) \in U$, an open set $ \hat S_x \subset \R^{n_x}$ with $\hat x \in \hat S_x$, and a unique function $p(y,\zeta): U \to \hat S_x$ such that  {$\hat x=p (\hat y, \hat \zeta)$} and $L(y, \zeta, p(y,\zeta))=0$   for all $(y,\zeta) \in U$. Moreover, $p$ is differentiable, and its Jacobian at $(\hat{y}, \hat{\zeta}, \hat{x})$ equals 
 \begin{align}
 -\left ( J^{\hat y, \hat \zeta, \hat x}_{x} \right)^{-1}     \left[\begin{array}{cc}
    J^{\hat y, \hat \zeta, \hat x}_{y}\\
     J^{\hat y, \hat \zeta, \hat x}_{\zeta}
  \end{array}\right].
 \end{align}
 \end{theorem}  
 {From here on we call the main condition mentioned in Theorem~\ref{thm:ift} \emph{the Implicit Decision Rules (IDR) condition} for brevity:
\begin{definition}[The Implicit Decision Rules (IDR) condition] \label{def:IDR}
 We say that the IDR condition holds  for $\hat y\hspace{-0.05cm} \in \hspace{-0.05cm} S_y$, $\hat \zeta  \hspace{-0.05cm} \in \hspace{-0.05cm} \Omega$ if there exists  $\hat{x} \hspace{-0.05cm} \in \hspace{-0.05cm} \R^{n_x}$ such that $L(\hat y, \hat \zeta, \hat{x}) \hspace{-0.05cm} = \hspace{-0.05cm} 0$ and the Jacobian $J^{\hat y, \hat \zeta, \hat x}_{x} $  is non-singular. We say that the IDR condition holds on sets $\hat{S}_y$, $\hat{\Omega}$ if it holds for each $y\hspace{-0.05cm} \in \hspace{-0.05cm} \hat{S}_y$, $\zeta \hspace{-0.05cm} \in \hspace{-0.05cm} \hat{\Omega}$.
\end{definition}}
 {We assume that it is possible to check quickly if the IDR condition holds for given $\hat y, \hat{\zeta}$ since one can solve equalities~\eqref{eq:L0} rather efficiently in practice. This paper was inspired by engineering problems, such as ACOPF. For these problems, Newton's methods~\cite{zimmerman2010matpower} are usually well-developed, fast, and precise enough to work for realistic problem instances.} Theorem~\ref{thm:ift} implies that, if the  IDR condition holds for some $(\hat y,\hat \zeta)$,  a unique local second-stage decision rule $x=p(y,\zeta)$ as in~\eqref{pr:0gen_poSubs} exists on some open set around $(\hat y, \hat{\zeta})$, and the first-order Taylor approximation of this rule is 
\begin{align}
    \hat p^{\hat y, \hat \zeta, \hat x} (y,\zeta) : = \hat x- \left( J^{\hat y, \hat \zeta, \hat x}_{x} \right)^{-1}\left ( J^{\hat y, \hat \zeta, \hat x}_{y}(y-\hat y) +J^{\hat y, \hat \zeta, \hat x}_{\zeta}(\zeta-\hat \zeta)\right ), \  {\text{ where } \hat x =p(\hat  y,\hat  \zeta)}. \label{eq:TaylImpl}
\end{align}
The above result applies to $\hat y \in S_y$ for which the Jacobian with respect to $x$ is non-singular. This condition is in general desirable in applications to dynamical systems, which are the main target of this paper. For such applications, solutions with singular Jacobians  may be physically unstable, see, for example, \cite{wang2000investigation}.

 {
The idea is to work on small subsets of $S_y$ and $\Omega$ and use piecewise affine decision rules based on the above Taylor approximations within subsets. Fix some subset of $S_y$ and $\hat{y}\in \hat S_y$. For $M_\zeta>0$,  partition  $\Omega$ into $\Omega^1,\dots, \Omega^{M_\zeta}$  and pick some $\hat{\zeta}^1 \in \Omega^1,\dots, \hat{\zeta}^{M_\zeta}\in \Omega^{M_\zeta}$. If $\hat S_y$ is small enough and  the IDR condition holds for $\hat{y}$ and $\hat \zeta^1,\dots,\hat \zeta^{ M_\zeta}$,  the implicit decision rules from Theorem~\ref{thm:ift} exist on $\hat S_y$ and some open subset of $\Omega$. Now, in problem~\eqref{pr:0gen_poSubs}, we impose the inequalities for each $\Omega^k$  substituting the implicit rules   with the corresponding approximation~\eqref{eq:TaylImpl}.  The procedure results in an almost standard ARO approximation with piecewise affine decision rules  (see problem~\eqref{pr:TaylorInit} in Algorithm~\ref{alg:Taylor}). The  differences are that (i) $y$ is restricted to a subset, and (ii) the decision rules are known, we do not need to optimize over them.  We can solve the approximation and obtain  some solution $y^*$.  As the last step, we could check that  $y^*$ is feasible for some selected values of the uncertainty that are desirable in practice, for instance, that it is nominally feasible. Then we can repeat the procedure considering other subsets of $S_y$ and choose $y^*$ with the best objective value. \\
Intuitively, if the subsets of $S_y$ and $\Omega$ are small enough, the above approach should be able to  discard solutions $\hat y \in S_y$ for which no decision rules from Theorem~\ref{thm:ift} exist  and select solutions for which these rules exist. We develop this idea in the next corollary.} 
\begin{corollary}\label{cor:yfeas}   Consider an $\hat y \hspace{-0.05cm}\in \hspace{-0.05cm}S_y$,  a finite cover $\Omega \hspace{-0.05cm}= \hspace{-0.05cm}\bigcup_{k=1}^{ {M_\zeta}} \Omega^k$ and $\hat \zeta^1 \in \Omega^1,\dots, \hat \zeta^{M_\zeta} \in \Omega^{M_\zeta}$. 
 \begin{enumerate}[label=(\alph*)]
 \item Let  {the IDR condition} hold for  $\hat y$ on $\hat \zeta^1,\dots, \hat \zeta^{M_\zeta}$, and assume that equalities~\eqref{eq:L0} have no solution for  $\hat y$ on  some $\hat{\Omega} \subseteq \Omega$.  If each $ \Omega^k, k \hspace{-0.05cm}= \hspace{-0.05cm}1,\dots,M_\zeta $ is contained in a ball of radius $\epsilon \hspace{-0.05cm} >0$, then  {$\hat{\Omega}$ cannot contain a ball of radius larger than $2\epsilon$}.  \label{cor:yfeas1} 
 \item  Let  {the IDR condition} hold for  $\hat y$ on $\Omega$.  Then there exist  {a closed ball $\hat{S}_y\subset \R^{n_y}$ with $y\in \hat S_y$} and a cover $\Omega \hspace{-0.05cm}= \hspace{-0.05cm}\bigcup_{k=1}^{ {M_\zeta}} \Omega^k$  by finitely many compact subsets such that the unique decision rules from Theorem~\ref{thm:ift} exist on  products of $\hat S_y$ and each $\Omega^{1},\dots,  \Omega^{ {M_\zeta}}$. \label{cor:yfeas2} 
 \end{enumerate}  
\end{corollary}
 {\begin{proof}
We begin by proving part~\ref{cor:yfeas1}. In each $ \Omega^k, \  k=1,\dots,M_\zeta$ there is at least one $\zeta$ for which the IDR condition holds and therefore equalities~\eqref{eq:L0} have solutions.  Hence $\hat \Omega$ does not contain fully one of $\Omega^k$. To cover a ball of radius up to $r$ with balls of radius $\epsilon$ in such a way that none of the smaller balls is fully inside the larger ball, we need to have $r\le 2\epsilon$. Thus we cannot have a violating subset that contains a ball of radius larger than $2\epsilon$. Next we prove part~\ref{cor:yfeas2}. By Theorem~\ref{thm:ift}, for every $\zeta \in \Omega$ there is an open ball $B\in \R^{n_y} \times \R^{n_\zeta}$ with radius $\epsilon$ such that $(\hat y, \zeta) \in B$ and decision rules $x=p(\hat y,\zeta)$  exist on $B$. Now, consider an open ball $B_y$ around $y$ with radius $\eps/2$ and an open ball $B_\zeta$ around $\zeta$ with the same radius. Notice that $B_y\times B_\zeta \subset B$. We can cover  $\Omega$ by the balls $B_\zeta$ for every $\zeta \in \Omega$, and by compactness of $\Omega$ there will be a finite subcover.  Taking a closure of each ball in this subcover and intersecting it with $\Omega$, we obtain a finite compact cover $\Omega \hspace{-0.05cm}= \hspace{-0.05cm}\bigcup_{k=1}^{ {M_\zeta}} \Omega^k$. 
 Now, each of $\Omega^k$ corresponds to an open subset of $\R^{n_y}$ that contains $\hat y$ for which the decision rules exist by construction. Since there are finitely many such subsets, their intersection results in an open subset of $\R^{n_y}$ that contains $\hat{y}$ and the required closed ball $\hat{S}_y$ around it. 
\end{proof}}

\subsection{Approximation algorithm for general problems} \label{subsec:approx_gen}
We formalize the ideas from the previous section  in Algorithm~\ref{alg:Taylor}, and the approximation guarantee is presented in Theorem \ref{thm:Taylor}.
\begin{algorithm}[!htbp] \label{alg:Taylor}
\caption{Piecewise affine  approximations of problem~\eqref{pr:0gen_po} using Taylor series}
\smallskip
\linespread{0.6}\selectfont{}
\SetKwInOut{Input}{Input} 
\Input{Problem~\eqref{pr:0gen_po}, $\hat y^i \in  S^i_y \subseteq S_y , i=1, \dots,  {M_y}$ for some $  {M_y}$,  $\hat \zeta^{k} \in \Omega^{k} \subseteq \Omega$ $k=1,\dots, M_\zeta$ for some $ {M_\zeta}$, where  $\Omega \hspace{-0.05cm}= \hspace{-0.05cm}\bigcup_{k=1}^{ {M_\zeta}} \Omega^k$. All sets are compact.}
\For{$i=1,\dots, {M_y}$}{
 Check that  {the IDR condition} holds for $\hat{y}^i$ and all $\hat \zeta^1,\dots,\hat \zeta^{M_\zeta}$, denote the corresponding solutions by $\hat x^{i,k}, \ k=1,\dots, M_\zeta$. \\
 \If{ {The IDR condition} fails for $\hat{y}^i$ and some of $\hat \zeta^1,\dots,\hat \zeta^{M_\zeta}$}{  \smallskip
set $z_i:= \infty, \ , y^{*,i}:=  {\emptyset}$}
\Else{ {Solve \textsc{Problem} \ref{pr:TaylorInit}:
(Piecewise affine approximation of problem~\eqref{pr:0gen_po} on $S^i_{y}$) \vspace{-0.3cm}
\begin{subequations}  \label{pr:TaylorInit} 
\begin{align}  
z_i = \inf_{y} \quad &   f(y) \nonumber  \\[-0.2cm]
 \quad \st \quad &  y \in S^i_{y}\nonumber   \\[-0.1cm]
& \hat p^{\hat y^{i}, \hat \zeta^{k}, \hat x^{i,k}}(y,\zeta) \in S_x \text{ for all }  \zeta \in \Omega^{k}, \ k\in [ {M_\zeta}] \label{ineq:ximpl} \\[-0.1cm]
& G_m \left (y,\zeta, \hat p^{\hat y^{i}, \hat \zeta^{k}, \hat x^{i,k}}(y,\zeta) \right )\ge 0  
 \text{ for all }  \zeta \in \Omega^{k} , \ m\in [m_{in}],  \ k  \in [ {M_\zeta}]  
 \label{ineq:Gimpl} \\[-0.7cm] \nonumber
\end{align} 
\end{subequations} } 
where $\hat p^{\hat y^{i}, \hat \zeta^{k}, \hat x^{i,k}}(y,\zeta)$ is defined  in \eqref{eq:TaylImpl}. Let $y^{*,i}$ be the optimal solution to~\eqref{pr:TaylorInit}. \smallskip

Save $z_i, \ y^{*,i}$ setting  $z_i:= \infty, \  y^{*,i}:=  {\emptyset}$ if problem~\eqref{pr:TaylorInit} is infeasible.}} \vspace{-0.1cm} 
Choose $i^*:= \text{arg} \min_{i=1}^{ {M_y}} z_i, \ z^*:= z_{i^*},  \ y^*:=y^{*,i^*}$. \\
  {For  $\hat \zeta^{k}, \ k=1,\dots, {M_\zeta}$, check if $y^*$ is feasible for problem~\eqref{pr:0gen_po} restricting $\Omega$ to $\hat \zeta^{k}$.}\\
  Return $z^*, \ y^*$ if the  check in step 9  is positive, return $z=\infty, \ y^* =  {\emptyset}$ otherwise. \smallskip
\end{algorithm} 
The next theorem shows that Algorithm~\ref{alg:Taylor} provides approximation guarantees and is able to find solutions robust in a ``strong" sense (such that the inequality constraints hold with a  margin). 
\begin{theorem}[Approximation guarantees of Algorithm~\ref{alg:Taylor} for problem~\eqref{pr:0gen_po}] \label{thm:Taylor}   Consider a $\delta>0$, and let $y^*  \neq \emptyset$ be a  solution returned by Algorithm~\ref{alg:Taylor}.
\begin{enumerate}[label=(\alph*)] 
\item[] \hspace{-1.4cm}  {\textbf{Guarantees limiting infeasibility} }
\item There exists $\epsilon>0$ such that, if  {$S^{i^*}$} and each subset in the partition  $\Omega \hspace{-0.05cm}= \hspace{-0.05cm}\bigcup_{k=1}^{ {M_\zeta}} \Omega^k$ fit in balls of radius $\epsilon$, then the following conditions hold on $\Omega$, except for possibly a subset  {where no ball of size larger than  $2\epsilon$ fits}: (i) equality constraints~\eqref{eq:L0} hold, and (ii) each inequality constraint~\eqref{ineq:G0},~\eqref{ineq:x0} is violated by at most $\delta$. \label{thm:Taylor1}
\item   {Let $y^*$ be not feasible for problem~\eqref{pr:0gen_po} on $\hat{\Omega} \subset \Omega$ whose Lebesgue measure constitutes at least $\gamma \hspace{-0.05cm}  < \hspace{-0.05cm}  1$ of that of $\Omega$.  If  $M_\zeta  \hspace{-0.05cm}  \ge  \hspace{-0.05cm}  \tfrac{\log{\delta}}{\log(1-\gamma)}$ and $\hat \zeta_1, \dots, \hat \zeta_{ {M_\zeta}}$ are sampled uniformly i.i.d. \hspace{-0.05cm}on $\Omega$,  the probability that  $y^*$ passes the check in step 9 of Algorithm~\ref{alg:Taylor} is not larger than $\delta$.\label{thm:Taylor_add} } 
\item[] \hspace{-1.4cm}  {\textbf{Guarantees for feasible solutions} }
\item  Consider the set $S^\delta_y\subseteq S_y$ of all solutions $y$ to problem~\eqref{pr:0gen_po} which (i) satisfy  {the IDR condition} on $\Omega$ and (ii) are such that for any $\zeta \in \Omega$ constraints~\eqref{ineq:G0},~\eqref{ineq:x0} hold strictly with a margin $\delta$ (i.e., the inequalities $``\ge \delta"$ hold). Then there is a  {subset} of $S_y$ and a finite cover of $\Omega$ for which Algorithm~\ref{alg:Taylor} will return an optimal $y^*\in S^\delta_y$.   \label{thm:Taylor2}
\end{enumerate}
 \end{theorem}
\begin{proof} 
First, consider statement~\ref{thm:Taylor1}. We know that step  {2} in the algorithm did not fail for $y^*$ and the corresponding partition $\Omega=\bigcup_{k=1}^{ {M_\zeta}} \Omega^k$.  If $\eps$ is small enough, then by Corollary~\ref{cor:yfeas}~\ref{cor:yfeas1} there are decision rules $x=p( {\hat y^{i^*}},\zeta)$   on $\Omega$, except for possibly a subset  {where no ball of size larger than $2\epsilon$ fits.  A small enough $\epsilon$ ensures also that $S^{i^*}_y$ is such that  the decision rules are valid for both  $\hat y^{i^*}$ and $y^*$.} Thus part (i) of statement~\ref{thm:Taylor1} follows. Now, consider inequality constraints~\eqref{ineq:G0}~\eqref{ineq:x0}.  By assumption~\eqref{as:basic}~\ref{as:shapeZeta}, these constraints are continuous, and there are finitely many of them. Since we work on compact sets, Taylor approximation of each decision rule $\hat p^{\hat y^{i^*}, \hat \zeta^k, \hat x^{i^*,k}}(y^*,\zeta)$~\eqref{eq:TaylImpl} is close enough to $p(y^*,\zeta)$ on $ {S^{i^*}_y}$ and $\Omega^k$ so that the approximated inequality constraints~\eqref{ineq:Gimpl} in problem~\eqref{pr:TaylorInit} are $\delta$-close to the original inequality constraints~\eqref{ineq:G0},~\eqref{ineq:x0}  in problem~\eqref{pr:0gen_po}. Thus  the original constraints cannot be violated  {(within the subsets of $\Omega$ where the rules exist)} by more than $\delta$, so $y^*$ possesses property (ii).  {To prove statement~\ref{thm:Taylor_add}, notice that with i.i.d. uniform sampling, the probability that all constraints of problem~\eqref{pr:0gen_po} are satisfied for one sample $\zeta$ is at most $1-\gamma$. Hence the probability that they are satisfied for $M_\zeta$ samples is at most $(1-\gamma)^{M_\zeta}$, and if $M_\zeta\ge \tfrac{\log{\epsilon}}{\log(1-\gamma)}$, then the latter probability is not larger than $\epsilon$.} 
Finally, consider statement~\ref{thm:Taylor2}. By Corollary~\ref{cor:yfeas}, there exist a compact subset of $S_y$ that contains $y^*$ and a compact finite cover of $\Omega$ where  the implicit decision rules used by Algorithm~\ref{alg:Taylor} are valid.  By the above argument for part (ii) of statement~\ref{thm:Taylor1}, if we make the subset of $S_y$  and each subset in the cover of $\Omega$ small enough, the Taylor approximations are so close to the actual decision rules that the inequality constraints cannot differ from their approximations by more than $\delta$, and thus Algorithm~\ref{alg:Taylor} will choose $y^*$  {if we restrict the search to subsets from $S_y^\delta$.}
\end{proof}
To our knowledge, Algorithm~\ref{alg:Taylor} is the first ARO algorithm which uses the second-stage policies defined by implicit functions, represented by equalities~\eqref{eq:L0}.  {Theorem~\ref{thm:Taylor}\ref{thm:Taylor_add} speaks about uniform sampling from $\Omega$, which is a numerically efficient procedure since $\Omega$ is an ellipsoid \cite{dezert2001efficient}.}

The main limitation of the algorithm in practice is the need to choose the subsets of $S_y$ and $\Omega$. This is not  trivial since inefficient partitioning may lead to long running times and numerically unstable problems. Constructing good subsets is a topic for separate research. For instance, to split $\Omega$, one can use the results from \cite{splittingUnc}.  In this paper, we are especially interested in testing  Algorithm~\ref{alg:Taylor} as a proof-of-concept for robust ACOPF. This problem has an additional structure:  control variables $y$ appear linearly and separately in $L$. In the next subsection, we introduce a version of  Algorithm~\ref{alg:Taylor} that uses this property  {and allows working on larger subsets of $S_y$}.

\subsection{Approximation algorithm for problems where $L$ is linear in the control variables} \label{subsec:approx_our}

In this subsection, we assume that the control variables $y$ appear linearly and separately in the equality constraints $L(y,\zeta, x)=0$ . That is, we can write
\begin{align}
  L(y,\zeta, x) =L_1(x,\zeta) +J_y y,   \label{eq:Liny}
\end{align}
where $L_1$ is a polynomial function and $J_y$ is a constant matrix. As a result, the Jacobian of $L(y,\zeta,x)$ with respect to $y$ is $J_y$, and other Jacobians do not depend on $y$. Denote the Jacobians with respect to $\zeta$ and $x$ for any $y$ and some $(\hat \zeta, \hat  x)$ by $J^{\hat \zeta, \hat x}_{\zeta}$ and $J^{\hat \zeta, \hat x}_{x}$, respectively.  Then the first-order Taylor approximation~\eqref{eq:TaylImpl} in $(\hat y, \hat \zeta, \hat x)$ where $L(\hat y, \hat \zeta, \hat x)=0$ simplifies to
\begin{align}
   \hat p^{\hat \zeta, \hat x}(y,\zeta):=  \ & \hat x- \left( J^{\hat \zeta, \hat x}_{x} \right)^{-1}\left ( J_{y}(y-\hat y) +J^{\hat \zeta, \hat x}_{\zeta}(\zeta-\hat \zeta)\right ) \nonumber \\
   = \ & \hat x- \left( J^{\hat \zeta, \hat x}_{x} \right)^{-1}\left ( J_{y}y + J^{\hat \zeta, \hat x}_{\zeta}(\zeta-\hat \zeta) +   {L}_1(\hat \zeta, \hat x) \right ). \label{eq:TaylImpl2}
\end{align}
Observe that  {the Jacobian $J^{\hat \zeta, \hat x}_{x}$ and} the above Taylor approximation are independent of $\hat y$. That is, the same approximation could be valid for various subsets of 
$S_y$ as long as the decision rules for these subsets and some $\Omega^k \subseteq \Omega $ lead to the same subset $\hat S_x \subset \R^{n_x}$. Hence it could be beneficial  to use subsets of $S_x$ instead of subsets of $S_y$,   
 reconsidering Algorithm~\ref{alg:Taylor} as presented in Algorithm~\ref{alg:main}.

\begin{algorithm}[!htbp] \label{alg:main}
\caption{{Piecewise affine approximations of problem~\eqref{pr:0gen_po} where $L$ is linear in $y$}}
\smallskip
\linespread{0.6}\selectfont{}
\SetKwInOut{Input}{Input} 
\Input{Problem~\eqref{pr:0gen_po},  {$\gamma>0$},  $\hat x^j \in S_x , j=1, \dots,  {M_x}$ for some $  {M_x}$,  $\hat \zeta^{k} \in \Omega^{k} \subseteq \Omega$ $k=1,\dots, M_\zeta$ for some $ {M_\zeta}$, where  $\Omega \hspace{-0.05cm}= \hspace{-0.05cm}\bigcup_{k=1}^{ {M_\zeta}} \Omega^k$. All sets are compact.}
Solve \textsc{Problem} \ref{pr:split} (Approximation of problem~\eqref{pr:0gen_po} where $L$ is linear in $y$)  \vspace{-0.2cm}
\begin{subequations}\label{pr:split} 
\linespread{0.6}\selectfont{}
\begin{align}
z^* = \inf_{y} \quad &   f(y) \nonumber  \\[-0.2cm] 
 \quad \st \quad &  y \in S_{y}\nonumber   \\
 & \text{\textbf{For} } k=1,\dots, {M_\zeta}:  \nonumber \\
& \text{At least one of conditions (j) for } j=1,\dots, {M_x} \text{ holds:} \nonumber \\[-0.1cm] 
 \text{} (j):& \begin{cases} 
 \text{ {The IDR condition} holds for some } \hat y \in S_y \text{ on } \hat \zeta^ {k}\text{ with solution } \hat x^ {j},  \\[-0cm]
  { \| \hat p^ {{\hat \zeta^k, \hat x^j}}(y,\hat \zeta^k) -\hat x^j) \| \le \gamma  \quad (**)} \\
 {  \hat p^ {{\hat \zeta^k, \hat x^j}}(y,\zeta) \in S_x  \text{ for all } \zeta \in \Omega^k } \\
 G_m \left (y,\zeta, \hat p^ {{\hat \zeta^k, \hat x^j}}(y,\zeta) \right )\ge 0  \text{ for all } m\in [m_{in}],  \zeta \in \Omega^k, \end{cases} \label{ineq:Gimpl2} \\[-0.8cm]  \nonumber 
\end{align} 
\end{subequations} 
where $\hat p^ {{\hat \zeta^k, \hat x^j}}(y,\zeta)$ is defined  in \eqref{eq:TaylImpl2}. Let $y^*$ be the optimal solution to problem~\eqref{pr:split}.\\
   {For  $\hat \zeta^{k}, \ k=1,\dots, {M_\zeta}$, check if $y^*$ is feasible for problem~\eqref{pr:0gen_po} restricting $\Omega$ to $\hat \zeta^{k}$.} \\
  Return $z^*, \ y^*$ if the  check in step 2  is positive, return $z=\infty, \ y^* = \emptyset$ otherwise. \smallskip
\end{algorithm}
The next Corollary~\ref{cor:main} shows that  {the performance of Algorithm~\ref{alg:main} is similar to the one of Algorithm \ref{alg:Taylor}}.
\begin{corollary}[Approximation guarantees of Algorithm~\ref{alg:main}] \label{cor:main}  Let $y^*  \neq \emptyset$ be a  solution returned by Algorithm~~\ref{alg:main}. 
\begin{enumerate}[label=(\alph*)] 
\item  {The result of Theorem~\ref{thm:Taylor}~\ref{thm:Taylor_add} holds for $y^*$, and there is $\eps>0$ such that if $\gamma\le \epsilon$} and $\Omega^1,\dots,\Omega^{M_\zeta}$ fit in balls of radius at most $\epsilon$, the result of Theorem~\ref{thm:Taylor}~\ref{thm:Taylor1} holds for $y^*$. \label{cor:main1}
\item Define $\eps=\max_{\zeta_1,\zeta_2\in \Omega} \|\zeta_1-\zeta_2\|$. If $\epsilon$ is small enough, then   {using $M_x=M_\zeta = 1$} (no need to partition $\Omega$)  {will provide the guarantees of   Corollary~\ref{cor:main}~\ref{cor:main1} if $y^*\neq \emptyset$.} \label{cor:main2}
\end{enumerate}
 \end{corollary}
\begin{proof}  {We begin by proving statement \ref{cor:main1}. The result of Theorem~\ref{thm:Taylor}~\ref{thm:Taylor_add} holds by step 2 of Algorithm~\ref{alg:main}, which coincides with step 9 of Algorithm~\ref{alg:Taylor}. Now, move to the result of Theorem~\ref{thm:Taylor}~\ref{thm:Taylor1}. 
Denote the set of all  solutions $y\in S_y$ to~\eqref{eq:L0} for some $\hat \zeta^k, \ k=1,\dots, M_\zeta$ and $\hat x^{j}, j=1,\dots, M_x$  by $S^{j,k}_y$.  In each condition (j) in problem~\eqref{pr:split}, we verify that the IDR condition holds for some $y\in S^{j,k}_y$ and $\hat \zeta^k$ resulting in the second-stage solution $\hat x^{j}$. Thus the IDR condition holds for all $y\in S^{j,k}_y$ on $\hat{\zeta}^k$ since the Jacobian $J^{\hat \zeta^k, \hat x^{j}}_x$  does not depend on $y$. If the verification is successful, a second-stage rule exists for an open ball around any $y\in S^{j,k}_y$ and $\hat \zeta^k$. Since $S^{j,k}_y$ is compact, it has a finite cover by those balls. Consider the union of the balls in this cover, which is an open set containing  $S^{j,k}_y$, and look at the complement of this union. The distance between this complement and $S^{j,k}_y$ is positive since they are disjoint, $S^{j,k}_y$ is compact and the complement is closed. Hence, there is some $\epsilon>0$ such that if $\gamma \le \epsilon$, then the restriction on the Taylor approximation in \eqref{ineq:Gimpl2} (**) pushes $y^*$ to be in the earlier constructed  open cover of $S^{j,k}_y$, so there exist second-stage rules for $y^*$ around $\hat \zeta^k$. Moreover, we can adjust $\gamma$ and and the size of $\Omega^k$ to make sure that the Taylor approximations of the inequality constraints are $\delta$-close to the original constraints. These are the conditions we used to prove Theorem~\ref{thm:Taylor}~\ref{thm:Taylor1}, thus the corresponding result holds for $y^*$. } 
 Finally,  statement \ref{cor:main2} follows directly from statement~\ref{cor:main1} under the following consideration. If $\Omega$ is small enough, then  statement \ref{cor:main1} applies to the whole $\Omega$. By construction, Algorithm~\ref{alg:main} needs to find at least one subset of $S_x$ where the feasible decision rules apply. Thus, if we restrict the algorithm to one $\hat x\in S_x$ and obtain $y^*\neq \emptyset$, the conditions of Corollary~\ref{cor:main}~\ref{cor:main1} are satisfied. 
\end{proof}

 {Algorithm~\ref{alg:main} could be  more efficient than Algorithm~\ref{alg:Taylor} since it could capture more values of $y\in S_y$, even if it looks around  only one solution $\hat{x}$.} However, we cannot guarantee that a strictly feasible solution $y^*$ as defined in Theorem~\ref{thm:Taylor}~\ref{thm:Taylor2} will be found under some conditions as we cannot control $y$ as precisely as in Algorithm~\ref{alg:Taylor}. Also, to say that \emph{at least} one of the  conditions (j) in~\eqref{ineq:Gimpl2} holds for each $\Omega^k$, one may have to introduce additional binary control variables.

By Corollary~\ref{cor:main}~\ref{cor:main2}, if $\Omega$ is small enough for the chosen approximation error, it suffices to find one good small subset of $S_x$, and there is no need to partition $\Omega$. The nominal problem~\eqref{pr:nom} for which we want to find a robust solution is usually feasible and, moreover, its optimal solution is usually ``stable" in the sense that the Jacobian in~\eqref{ineq:Gimpl2} is non-singular. Thus, for small uncertainty sets, one could start from some nominal feasible control solution $\hat y^1$ at hand (e.g.,  the original optimal solution), find the corresponding state variables solution  $\hat x^1$ such that $L(\hat y^1, 0, \hat x^1)=0$ and restrict Algorithm~\ref{alg:main} to work around $\hat x^1$ as in (**) in~\eqref{ineq:Gimpl2}. If the algorithm finds a solution $\hat y^2 \neq \hat y^1$, one could set $\hat x^2$ such that $L(\hat y^2,0,\hat x^2)=0$ and repeat the procedure trying to find a better solution.  {Notice that a similar procedure would be valid for Algorithm~\ref{alg:Taylor} too, by Theorem~\ref{thm:Taylor}~\ref{thm:Taylor1}, but we would construct subsets around $y^1, y^2$ in $S_y$ instead of subsets around $x^1, x^2$   in $S_x$.  }

 {Now, we have constructed two algorithms, one is suitable for  general problems~\eqref{pr:0gen_po} (Algorithm~\ref{alg:Taylor}), and one is tailored to problems~\eqref{pr:0gen_po} linear in control variables (Algorithm~\ref{alg:main}). To solve the optimization problems~\eqref{pr:TaylorInit} and~\eqref{pr:split} in the algorithms, we need to deal with quadratic inequality constraints under ellipsoidal uncertainty, which is the topic of the next section. }

\section{Quadratic inequality constraints under uncertainty} \label{sec:Ineq} 
In this section, we  {solve problem~\eqref{pr:TaylorInit}  in Algorithm~\ref{alg:Taylor} or problem~\eqref{pr:split} in Algorithm~\ref{alg:main}, depending on which algorithm is used.} 

Our first step is to eliminate  $\zeta$ from~\eqref{ineq:Gimpl},~\eqref{ineq:ximpl} and~\eqref{ineq:Gimpl2}. The constraints in question have the following general form:
\begin{align}
h(y,\zeta) \ge 0 \text{ for all } \zeta \in \Omega. \label{ineq:gen}
\end{align}  
For instance, we can have $h(y,\zeta) :=G_i \left (y,\zeta,\hat p^{\hat y, \hat \zeta, \hat x}(y,\zeta)\right )$ for some $ i\in [m_{in}]$. 
Since~\eqref{ineq:gen} concerns the non-negativity of a polynomial on the set $\Omega$, whose description does not involve $y$, in the next subsection~\ref{subsec:uncert} we treat $y$ as a parameter and answer the question of when a polynomial of degree at most two in $\zeta$ is non-negative on $\Omega$. The answer leads to a reformulation of \eqref{ineq:gen} that contains no $\zeta$. In this way, we eliminate the uncertainty $\zeta$ from the problems in question. After that, we stop treating $y$ as a parameter, and $y$ becomes the only variable left in the reformulated problems ~\eqref{pr:TaylorInit} and~\eqref{pr:split}, which we solve in subsection~\ref{subsec: proj}. 

The rest of this section is needed if the initial problem~\eqref{pr:0gen_po} contains \emph{inequalities} that are non-linear in $y$ or $\zeta$. If it turns out that  {problem~\eqref{pr:TaylorInit}  in Algorithm~\ref{alg:Taylor} or problem~\eqref{pr:split} in Algorithm~\ref{alg:main}} is linear in $y$ and $\zeta$, then classical robust optimization techniques apply to the corresponding problem. In this case one can use a standard reformulation of a linear constraint under the ellipsoidal uncertainty  (see,e.g., \cite{BenTal2004}) and rewrite the problem as a second-order cone program.  {Finally, from here on, we focus on problem~\eqref{pr:split} for brevity, but the method is valid for problem~\eqref{pr:TaylorInit} as well.}

\subsection{Eliminating the uncertainty from the problem} \label{subsec:uncert}  
  Under Assumption~\ref{as:basic},  all inequalities in problem~\eqref{pr:0gen_po} are of degree at most two  in $x$ and $\zeta$, hence so is $h$ in~\eqref{ineq:gen}: 
\begin{align}
h (y,\zeta):=\zeta\tr A  \zeta+ (y\tr B +b \tr ) \zeta + c  \tr y +d  +g (y) \ge 0, \label{eq:quadG}
\end{align}
for some given parameters $A , \ B , \ b , \ c , \ d $, and a function $g $ which contains all monomials that are non-linear in $y$. The precise values of the parameters and the form of $g$ depend on the functions in~\eqref{ineq:Gimpl2} and can be obtained directly from those functions.  Recall that $\Omega$ has the form~\eqref{def:OmegaQuad}. Such a combination of $h$ and $\Omega$ allows reformulating \eqref{ineq:gen} and eliminating the uncertainty from it due to the \emph{S-lemma}.
\begin{proposition}[\cite{Yakubovich}] \label{prop:slemma}
Let $h,g\in \R_2[{\zeta}]$ and suppose there is ${\zeta} \in \R^{n_{{\zeta}}}$ such that $g({\zeta})<0$. Then the following two statements are equivalent:
\begin{enumerate}
\item $h({\zeta})\ge0$ for all $\zeta\in \R^{n_{{\zeta}}}$ such that $g({\zeta})\ge 0$.
\item There is $\lambda\in \R_+$ such that $h(\zeta)-\lambda g({\zeta}) \ge 0$ for all ${\zeta}\in \R^{n_{{\zeta}}}$.
\end{enumerate} 
\end{proposition} 
S-lemma is well known in robust optimization (see~\cite{BeTal2009}) but is usually applied to convex problems in $y, \zeta$ while we use it for a general quadratic constraint~\eqref{ineq:gen}.

\begin{proposition} \label{prop:slem}
Constraint~\eqref{ineq:gen} with $h$ as in~\eqref{eq:quadG} and the uncertainty set~\eqref{def:OmegaQuad} holds if and only if there exist $\lambda , \gamma $ such that 

\begingroup
\begin{subequations}
\linespread{0.6}\selectfont{}
\begin{align*}
  &  \begin{bmatrix}
 \gamma  +c  \tr  y+d +\lambda  r &  \tfrac{1}{2}(y\tr B +b \tr+\lambda  \sigma \tr )\\
\tfrac{1}{2}( B  \tr y+b +\lambda  \sigma) &   \lambda  \Sigma+ A 
\end{bmatrix}   \succeq 0\\
& \lambda  \ge 0,  \\
& g (y)=\gamma  ,
\end{align*}
\end{subequations}
\addtocounter{equation}{-1}
\endgroup
where $\Sigma, \sigma$ and $r$ are the  parameters from~\eqref{def:OmegaQuad}, and other parameters are defined in~\eqref{eq:quadG}. 
\end{proposition}
\begin{proof}
We have 
\begin{align*}
 & h (y,\zeta) \ge 0 \text{ for all } \zeta \in \Omega \\
\overset{\text{\tiny{Prop.~\ref{prop:slemma}}}}{\iff}& \zeta\tr A  \zeta+ (y\tr B +b \tr ) \zeta + c  \tr y +d  +g (y)  \ge 0 \text{ for all } \zeta\in \R^{n_\zeta}: \zeta\tr \Sigma \zeta +\sigma \tr \zeta + r \le 0\\
\overset{}{\iff} & \zeta \tr A  \zeta+ (y\tr B +b \tr ) \zeta + c  \tr  y +d  +g (y) - \lambda  (-r- \zeta\tr \Sigma \zeta - \sigma \tr \zeta) = \   \left [ \begin{smallmatrix}
1\\ \zeta
\end{smallmatrix} \right ]  \tr S  \left [ \begin{smallmatrix}
1\\ \zeta
\end{smallmatrix} \right ], \\
& \quad  S \succeq 0, \ \lambda  \ge 0\\
\iff & \begin{bmatrix}
 \gamma  +c  \tr  y+d +\lambda  r &  \tfrac{1}{2}(y\tr B +b \tr+\lambda  \sigma \tr )\\
\tfrac{1}{2}( B  \tr y+b +\lambda  \sigma ) &   \lambda  \Sigma+ A 
\end{bmatrix}   \succeq 0, \ g (y)=\gamma , \ \lambda  \ge 0. & 
\end{align*}  
\addtocounter{equation}{-1}
\end{proof}
If some  inequalities in problem~\eqref{pr:0gen_po} are linear in $x$ and $\zeta$, then $h$ in~\eqref{ineq:gen} is linear in $\zeta$ as well. In this case, the constraint considered in this section simplifies to the classical linear robust constraint. It is well known how to eliminate uncertainty from such constraints; see, e.g., \cite{BeTal2009}. Each linear constraint under the ellipsoidal uncertainty in~\eqref{def:OmegaQuad} will become a second-order cone constraint after eliminating the uncertainty.

\subsection{Alternating projections algorithm} \label{subsec: proj} 
In Section~\ref{subsec:uncert}, we eliminate $\zeta$ from  constraints \eqref{ineq:Gimpl2} in problem~\eqref{pr:split} and add auxiliary control variables (e.g., $\lambda, \gamma$ in Proposition~\ref{prop:slem}). { By eliminating $\zeta$, we {reformulate} problem~\eqref{pr:split}. This reformulation can be written} in a general way as
\begingroup
\let\theequation\theproblem
\begin{problem}
\linespread{0.6}\selectfont{}
\begin{subequations} \label{pr:reform_split}
\begin{align}
z_j =\inf_{\substack{y}} \quad &   f(y) \label{obj:reform_split1}\\
\st \quad & H(y)=\gamma \label{const:reform_split2} \\
& F \left[ \begin{smallmatrix} y\\ \gamma \end{smallmatrix} \right ]  \in \cC, 
 \label{const:reform_split3}   
\end{align}
\end{subequations}\end{problem}\endgroup
\noindent where $\gamma$ are auxiliary control variables as in Proposition~\ref{prop:slem}, $y$ includes other original and auxiliary control variables, $H$ is a polynomial mapping, $F$ is a linear mapping, and $\cC$ is a proper semialgebraic cone amenable for optimization (e.g., positive semidefinite cone, second-order cone, or the non-negative orthant). 

Under Assumption~\ref{as:basic}, the equality constraints are quadratic, the objective function is convex with degree at most two, and $\cC$ is the positive semideninite cone. The approach we use next is applicable for general polynomial equality constraints and semialgebraic cones $\cC$, therefore we introduce the formulation above. 
 
Problem~\eqref{pr:reform_split} would be a classical conic problem for which many solvers exist if we did not have those non-linear equality constraints. Hence, we first relax those complicating equalities and obtain a solution feasible for~\eqref{const:reform_split3}. Then we iteratively transform it into a solution feasible for the whole problem~\eqref{pr:reform_split}. Notice that problem~\eqref{pr:reform_split} cannot be unbounded as the feasible set of its nominal problem is compact by Assumption~\ref{as:basic}. We begin with a lower bound for the problem, which can be obtained using any relaxation of the polynomial equality constraints. One can use lifting techniques, e.g., SDP relaxations or the reformulation-linearization technique as in~\cite{sherali1992global}. The lower-bound relaxation might provide a solution that is infeasible for~\eqref{const:reform_split2}. To obtain a feasible solution from the solution to the relaxation, we use the alternating projection method as presented in Algorithm \ref{alg:quadAlg}. Let $z^U$ be an upper bound on $f(y)$ in problem~\eqref{pr:reform_split}. 
Define two sets:
\begin{align}
\linespread{0.6}\selectfont{}
    \A:= \ & \{(y,\gamma): \eqref{const:reform_split3} \text{ holds}, f(y)\le z^U \}, \label{set:A}\\
    \cB:= \ & \{(y,\gamma): \eqref{const:reform_split2} \text{ holds} \}. \label{set:B}
\end{align}
We denote $y:=(y^{nc}, y^{c})$, where $y^{nc}, y^{c}$ are the subsets of variables $y$ that are involved (resp. not involved) in the non-convex constraints~\eqref{const:reform_split2}. We use two subsets of the variables to speed up the algorithm by skipping the iterations in which only $y^c$ changes. \medskip

 \begin{theorem} \label{thm:quadAlg} Let $\A\cap \cB \neq \emptyset$. Then the following holds:
\begin{enumerate}[label=(\alph*)]
\item Algorithm~\ref{alg:quadAlg} stops after finitely many iterations. The algorithm can either find a feasible solution, or report an infeasible problem, or not be able to find a solution in the given number of iterations.  \label{thm:quadAlg0} 
    \item If  Algorithm~\ref{alg:quadAlg} starts at a point $(y_0,\gamma_0)$ that is sufficiently close to $\A\cap \cB$, then it stops in a point $(y^*,\gamma^*)$ that is $tol$-close to $ \A\cap \cB$ for $N$ large enough. \label{thm:quadAlg1}
    \item  If $\A$ has a non-empty interior and  $(y_0,\gamma_0)$ is sufficiently close to $\A\cap \cB$, Algorithm~\ref{alg:quadAlg} converges to a point in $ \A\cap \cB$ linearly. We say that a sequence $(a_k)_{k=1}^\infty$ converges \emph{linearly} to $a$ if $\lim_{k\to \infty}\tfrac{\|a_{k+1}-a\|}{\|a_{k}-a\|}<1$.  \label{thm:quadAlg2} 
    \item For any $\delta>0$ there exists $N$ large enough and $tol$, $(\nu_i)_{i=1}^N$ small enough such that the solution obtained by Algorithm~\eqref{alg:quadAlg} is $\delta$-close to a local minimum. \label{thm:quadAlg3} 
\end{enumerate}
\end{theorem}
\begin{proof}  
Item~\ref{thm:quadAlg0} follows by construction of the algorithm since we limit the number of iterations.    Item~\ref{thm:quadAlg1} follows from Assumption~\ref{as:basic} and Theorem 7.3 in~\cite{drusvyatskiyAP1} as the sets $\A$ and $\cB$ are closed and semialgebraic. Item~\ref{thm:quadAlg2} follows from the fact that $\A$ is convex. Hence the normal cone at any point in the interior of $A$ equals $\{0\}$, and thus $\A$ and $ \cB$ satisfy the conditions of Theorem~2.1 in~\cite{drusvyatskiyAP1}, which implies linear convergence. 
Next we prove item~\ref{thm:quadAlg3}. For a given $\nu$, let $S_\nu:=\A\cap \cB \cap \{(y,\gamma): f(y)\le \nu f(y^*)+(1-\nu)z^L \}$. Let $(y^i,\gamma^i)$ be the best solution of the algorithm at iteration $i$, and assume that it is not the closest local optimum denoted by $(\hat y,\hat \gamma )$. Then $f(y^i) < f(\hat y)$ and there exist $\nu_i>0$ such that  $S_{\nu_i} \neq \emptyset$ and $(y^i,\gamma^i)$ is close enough to $S_{\nu_i}$ to satisfy the conditions of item~\ref{thm:quadAlg1}.  Notice that the latter argument might fail if we require $(\hat y,\hat \gamma )$ to be a global optimum. Using $\nu_i$, the algorithm finds the next point $(y^{k+1},\gamma^{k+1})$ that is $tol$-close to $S_{\nu_i}$ and such that $f(y^i) < f(y^{i+1}) \le f(\hat y)$. Since the interval $(f(y^i),f(\hat y))$ is bounded, for any $\delta>0$ there are $N$ large enough, $tol$ small enough and a sequence $(\nu_i)_{i=1}^N$ such that $f(y^N)$  is $\delta$-close to $f(\hat y)$.
\end{proof}

\begin{algorithm}[] \label{alg:quadAlg}
\caption{Alternating projections algorithm for problem~\eqref{pr:reform_split}}
\medskip
\linespread{0.6}\selectfont{}
\SetKwInOut{Input}{Input} 
\Input{small $tol>0$, upper bound $z^U \ge z_j$, $(\nu_i)_{i=1}^N \in (0,1]$, $N\ge 1$} 
 Solve a lower bound relaxation of problem~\eqref{pr:reform_split} and denote its solution by $(y_0,\gamma_0)$ and its objective value by $z^L$\;
 \If{the lower bound relaxation is infeasible,}{Lower-bound relaxation is infeasible, stop, return ``Problem~\eqref{pr:reform_split} infeasible"}
 Set $i=0$ and choose $(y_1,\gamma_1)$ such that $\|(y^{nc}_0,\gamma_0)-(y^{nc}_1,\gamma_1)\| > tol$\;
\While{$\|(y^{nc}_0,\gamma_0)-(y^{nc}_1,\gamma_1)\| > tol$ and $i \le  N$}{
Set $i:=i+1$\; 
Project on $\cB$: 
  $(y_1,\gamma_1):=  (y_0,H(y_0))$\;
  Project on $\A$: $ (y_0,\gamma_0): = \text{arg} \min_{(y,\gamma)\in \A} \|(y^{nc},\gamma)-(y^{nc}_1,\gamma_1)\|  $\;
\If{ $\|(y^{nc}_0,\gamma_0)-(y^{nc}_1,\gamma_1)\| < tol$,}{Find the best $y^{c}$ given $(y^{nc}_0,\gamma_0)$ by solving $y^{c}_0 = \text{arg} \min_{(y^{nc}_0,y^{c},\gamma_0)\in \A} f\left (y^{nc}_0,y^{c}\right )$\; 
Save $(y^*,\gamma^*):=(y^{nc}_0,y^{c}_0,\gamma_0)$ as the current best feasible solution\;
Try to find a feasible solution with a better objective value:\\ Update $z^U:=\nu_i f(y^*)+(1-\nu_{i})z^L \ \ $ (to decrease the upper bound)\; 
Update $\A:= \A \cap  \{(y,\gamma): f(y)\le z^U \} \ \ $ (to work with lower objective values)\;
Adjust $(y_1,\gamma_1):=\tfrac{1}{tol}(y_0,\gamma_0) \ \ $ (to proceed with the while loop)}
}
 \eIf{No feasible solution is obtained,}{ Problem~\eqref{pr:reform_split} could be infeasible, return ``Inconclusive, out of iterations"}{Return ``The best obtained solution is $(y^*,\gamma^*)$"}
\end{algorithm} 

In the next subsection, we demonstrate the performance of the combination of Algorithm~\ref{alg:main} with Algorithm~\ref{alg:quadAlg} to solve problem~\eqref{pr:0gen_po} on the ACOPF.  The experiments show that in the majority of  tested  cases  this combination of algorithms finds a robustly feasible solution   within 15 minutes of computational time for small uncertainty sets.

\section{Adjustable ACOPF with uncertain renewable generation and load demands} \label{subsec:acopf}

Optimal power flow (OPF) is one of the key optimization problems relevant to the operation of electric power systems. OPF solutions provide minimum cost operating points that satisfy both equality constraints termed the ``power flow equations'' which model the power system network and inequality constraints that impose limits on line flows, generator outputs, voltage magnitudes, etc. Accurately modeling the steady-state behavior of power systems requires the non-linear \emph{AC power flow equations}, which can be formulated as a system of quadratic equality constraints. 

Compounding the difficulties posed by the power flow non-linearities, rapidly increasing quantities of wind and solar generation are introducing significant amounts of power injection uncertainties into electric grids. To address these uncertainties, researchers have studied a wide range of stochastic and robust OPF problems~\cite{pscc2022survey}, many of which use the DC power flow approximation; see, e.g.,~\cite{vrakopoulou2013,ChanceBienstock}  
for several relevant examples. This linear power flow representation permits the application of stochastic and robust optimization techniques developed for linear programs. Alternative approaches replace the AC power flow equations with other more sophisticated approximations, such as the work in \cite{chaos1} and \cite{roald2018}, 
or convex relaxations, such as the work in \cite{venzke2017}. 
Such approaches can provide useful solutions in many contexts, particularly when the approximations are iteratively updated or adaptively adjusted. However, the quality guarantees from these approaches are provided with respect to the approximation or convex relaxation as opposed to the original non-convex ACOPF problem. Since power flow approximations and relaxations can introduce substantial errors relative to the non-linear AC power flow equations~\cite{purchala2005,cdc2016,baker2020}, the resulting solutions may lead to unacceptable constraint violations during operation in the physical system.

The power systems literature also includes approaches that directly address the non-linear AC power flow equations. These approaches can provide high-quality solutions in certain instances but may be limited to special classes of problems, such as systems that satisfy restrictive requirements on the power injections at each bus as in \cite{RobustACOPF}. Other approaches use scenario-based techniques that enforce feasibility for selected uncertainty realizations, 
possibly obtained via subproblems that compute worst-case uncertainty realizations with local solvers as in~\cite{capitanescu2012} or convex relaxations as in~\cite{AdaptRobACOPF}. Certifying robustness with such approaches is challenging due to the possibilities of local solutions and inexact relaxations. Rather than seeking the worst-case uncertainty realizations, the approach in~\cite{ACOPFRobFeas} instead bounds the worst-case impacts of the uncertainties with respect to each constrained quantity. While this approach provides guarantees regarding the satisfaction on the engineering inequality constraints, each iteration requires the solution of many computationally expensive subproblems. We also note recent work in~\cite{lee2020robust} that uses so-called ``convex restriction'' techniques (see~\cite{ConvexRestr1}) to compute robustly feasible ACOPF solutions. While promising, this approach is undergoing continuing development and requires specialization to the particular non-linearities in each class of problems.

Accordingly, new computational methods such as the framework proposed in this paper are needed to address power system optimization problems that model both uncertain power injections from renewable generators and non-linearities from the AC power flow equations. The proposed framework would most naturally fit into a system operator's generation scheduling processes as part of near-real-time (five-minute to hourly) generator setpoint computations via optimal power flow algorithms, which is a key application at the heart of power system operations. Since optimization in this case is repeated frequently during the day, the uncertain deviations in loads or generation are likely to be of small size. However, they can lead to infeasibility of the nominal solutions, thus an approach is needed that takes the uncertainty into account. This situation corresponds well with the result in Corollary~\ref{cor:main}~\ref{cor:main2} about small uncertainty sets. Extensions to incorporate binary variables would facilitate other applications, such as solving unit commitment problems for day-ahead scheduling that models generators' start up and shut down decisions. Thus, developing algorithms suitable for such extensions is an important direction for future work.

\subsection{Robust ACOPF formulation} \label{subsec:acopf_form}
Exploiting the polynomial representation of the AC power flow equations, we next apply the approach described in this paper to the robust ACOPF problem, beginning with our notation and the problem formulation. Consider a power network with the set of buses $N=\{1,\dots,n\}$ and the set of lines connecting these buses $E$. We denote the set of buses with generators by $G$ and the active and reactive power demand (load) at each bus $k\in N$ by $P^d_k$ and $Q^d_k$, respectively and denote the index of the reference bus by $s$. To implement thermal restrictions on the transmission lines, we impose line current limits, see~\cite{zimmerman2010matpower}. Our objective is to minimize the cost of power generation, which is one of classical objectives in OPF problems.  
Denote the active and reactive power injections due to load or generation fluctuation by $P_k^r$ and $Q_k^r$, respectively, for all $k\in N$. In the nominal ACOPF  without uncertainty, $P_k^r$ and $Q_k^r$ are known and fixed. Next we define the ACOPF problem as a quadratic optimization problem.
\begingroup\begin{problem} \label{pr:acopf}
\linespread{0.6}\selectfont{}
\begin{subequations}
\begin{align*}
z^{nom}=\inf_{x, P^g,Q^g} \quad & \sum_{k\in G} c_k^2 (P_k^g)^2+c_k^1 P_k^g+c_k^0 \\
\st \quad & P_k^{\min}\le P_k^g \le P_k^{\max} & \text{for all } k\in N \\
& Q_k^{\min}\le Q_k^g \le Q_k^{\max} & \text{for all } k\in N \\
& (V_k^{\min})^2\le x_k^2+x_{k+n}^2\le (V_k^{\max})^2 & \text{for all } k\in N \\
& \tra(Y_{lm}xx\tr)  \le S_{lm}^{\max} & \text{for all } \{lm\}\in E\\
& P_k^g+P^r_k= P^d_k +\tra(Y_kxx\tr) & \text{for all } k\in N \\
& Q_k^g+Q^r_k= Q^d_k+\tra(\bar Y_kxx\tr) & \text{for all } k\in N \\
& x(s+n)=0,
\end{align*}%
\end{subequations}%
\end{problem}
\addtocounter{equation}{-1}
\endgroup
\noindent where the last constraint sets the phase angle of the reference bus to zero. Now, let the active power fluctuations for each $k\in N$ be $P_k^r= \bar P_k^r +\zeta_k,$ where $\zeta$ represents the uncertainty. We assume that the power injection uncertainties from the load and generation at each bus $k\in N$ are modeled via a constant power factor $\cos \phi_k$ so that the reactive power fluctuations are
\begin{align*}
    Q_k^r=\bar Q_k^r+ \gamma_k \zeta_k, \  \gamma_k:= \begin{cases} 0& \text{ if } P^d_k=0,\\ \tfrac{\sqrt{1-\cos^2 \phi_k}}{\cos \phi_k} & \text{ otherwise}. \end{cases} 
\end{align*} 
Without loss of generality, we let $\bar P_k^r=\bar Q_k^r=0$, otherwise one can adjust the loads $P_k^d$ and $Q_k^d$. We denote by $\delta$ the total change in the active power losses due to the redistribution of power flows from the uncertain power injection fluctuations relative to the losses from the nominal power flows. Note that $\delta$ is typically near zero, as the losses themselves are usually small and the changes in losses are even smaller. For an operating point to be robustly feasible, the generators must account for the total change in the active power injections, $\sum_{i=1}^n \zeta_i-\delta$, associated with each uncertainty realization without leading to constraint violations. We adopt a ``participation factor'' model where this change in power injections is distributed among all generators according to a linear recourse policy with specified participation factors $\alpha_k$ for each generator $k$. Thus, for each $k\in N$, the actual active power generation consists of the nominal power $P^g_k$ and an adjustment in generation due to the uncertainty:
\begin{align}
& P_k^g- \alpha_k \left  (\sum_{i=1}^n \zeta_i -\delta \right  ),  \ \alpha_i\ge 0, \ \sum_{i=1}^n \alpha_i =1, \label{eq:policy}
\end{align}
Thus, when introducing uncertainty to problem~\eqref{pr:acopf}, we replace  $P_k^g$ in this problem with~\eqref{eq:policy}. We note that this model represents the steady-state behavior of widely used automatic generation control (AGC) (see \cite{jaleeli1992}) and is adopted in many robust and stochastic OPF formulations, e.g., those used by \cite{venzke2017}, \cite{roald2018}, and \cite{ACOPFRobFeas}. To model the uncertainty, we let  the uncertain parameters $\zeta$ belong to the region $\Omega=\{\zeta\in \R^{n_\zeta}: \zeta \tr \Sigma \zeta \le 1 \},$ where $\Sigma$ is a covariance matrix. That is, our uncertainty region is an ellipsoid centered on the point with no fluctuations.

For $k \in N$, we denote by $V^g_k:=x_k^2+x_{k+n}^2$ the squared voltage magnitude at bus $k$. Following traditional power system modeling practices, we consider three types of buses: PV, PQ and the reference bus. If $k$ is a PV bus, the active power $P^g_k$ and squared voltage magnitudes $V^g_k$ are set by the operator while the reactive power $Q^g_k$ can change. If $k$ is a PQ bus, then the active power and reactive power are fixed to constant values while the voltage magnitude can change. Without loss of generality, we assume that active and reactive power generation at PQ buses is zero, otherwise the loads can be adjusted. Finally, the operator selects the voltage magnitude at the reference bus while the active and reactive powers are free to vary. We also introduce a variable $t$ that denotes the worst-case upper bound on the active power on the reference bus. We use this  bound to estimate the worst-case objective value over the uncertain power injection fluctuations, as is typical in robust optimization problems. As a result, the control variables $y$ in the problem include $t$, $P^g_k$, where $k$ belongs to the set of PV buses, and $V^g_k$, where $k$ belongs to the union of PV and reference buses. 

Now we define the problem in the same form as \eqref{pr:0gen_po} to more easily use the results from the earlier sections. This yields the following:
\begin{align}
f(t,P^g,V^g)= & \sum_{k\in G \setminus \{s\}} c_k^2 (P_k^g)^2+c_k^1 P_k^g+c_k^0+ c_s^2 t^2+c_s^1 t+c_s^0, \label{def:f} \\
S_{y}= & \big \{(P^g,V^g):  P_k^{\min}\le P^g_k \le P_k^{\max} \  \text{ for all }   k \in G\setminus \{s\}, \nonumber \\
& \quad (V_k^{\min})^2\le V^g_k \le (V_k^{\max})^2  \  \text{ for all } k \in G, \ P^{\min}_s\le t \le P^{\max}_s \big \}, \label{def:Sy}\\
L(P^g,V^g,\zeta,x)=& \begin{bmatrix}
 P^d_k +\tra(Y_kxx\tr)-\zeta_k -P_k^g+ \alpha_k \sum_{i=1}^n \zeta_i & \text{for all } k\in G \setminus \{s\}\\
P^d_k +\tra(Y_kxx\tr)-\zeta_k+ \alpha_k \sum_{i=1}^n \zeta_i & \text{for all } k\in N \setminus G \\
 Q^d_k+\tra(\bar Y_kxx\tr)-\gamma_k\zeta_k &  \text{for all } k\in N \setminus G\\
 x_k^2+x_{k+n}^2-V^g_k & \text{for all } k\in G\\
 x(s)
\end{bmatrix}, \label{def:L} 
\end{align}
\begin{align}
G(t,P^g,V^g,\zeta,x)=& \begin{bmatrix}
 -P_s^{\min}+ P^d_s +\tra(Y_sxx\tr) -\zeta_s + \alpha_k \sum_{i=1}^n \zeta_i\\
 t-P^d_s -\tra(Y_sxx\tr) +\zeta_s- \alpha_k \sum_{i=1}^n \zeta_i   \\
 -Q_k^{\min}+Q^d_k+\tra(\bar Y_kxx\tr) -\gamma_k\zeta_k  & \text{for all } k\in G\\
 Q_k^{\max}-Q^d_k-\tra(\bar Y_kxx\tr)+\gamma_k \zeta_k  & \text{for all } k\in G\\
 \tra(Y_{lm}xx\tr) \le S_{lm}^{\max} & \text{for all } \{lm\}\in E\\
(V_k^{\max})^2- x_k^2+x_{k+n}^2 & \text{for all } k \in N \setminus G \\
x_k^2+x_{k+n}^2-(V_k^{\min})^2  & \text{for all } k \in N \setminus G 
\end{bmatrix}, \label{def:G}
\end{align}
\begin{align}
S_{x}= & \left \{x \in \R^{2n}: (V_k^{\min})^2\le x_k^2+x_{k+n}^2 \le (V_k^{\max})^2 \ \text{for all } k \in N \setminus G \right \} \label{def:Sx} 
\end{align}
In the next subsection, we run numerical experiments solving problem~\eqref{pr:0gen_po} with the inputs defined above and the instances from M{\sc atpower}~\cite{zimmerman2010matpower}. 

\subsection{Numerical results} \label{subsec:acopfNum}

In this section, we implement Algorithm~\ref{alg:main} and compare it with several algorithms used in the literature. All computations are done using MATLAB R2021a and Yalmip (see~\cite{Lofberg2004}) on a computer with the processor Intel\textsuperscript{\textregistered} Core\textsuperscript{\textregistered} i7-8665U CPU @ 1.90GHz and $16$ GB of RAM. Semidefinite programs are solved with MOSEK, Version 9.3.21 ~\cite{mosek}.

In the experiments, we need to choose the vector $\alpha$ and the parameters defining the ellipsoid $\Omega$. We assume that all uncertainty balances out on the reference bus since this is the default setup of  M{\sc atpower}. Therefore $\alpha_s=1$ and $\alpha_k=0$ for $k\neq s$.  \ \\
We implement Algorithm~\ref{alg:main} with the following inputs:
\begin{itemize}
    \item We do not split the uncertainty set $\Omega$ and set $\hat{\zeta} = 0$.
    \item   We consider only one $\hat x$: the initial feasible solution to the solution of the nominal problem~\eqref{pr:acopf} provided by M{\sc atpower}~\cite{zimmerman2010matpower}. For constraint (**) 
 in Algorithm~\ref{alg:main}, we set $\gamma=\sqrt{\sfrac{\|x_{j-1}\|}{10}}$ for $|N|<30$, and $\gamma=\sqrt{\sfrac{\|x_{j-1}\|}{30}}$ for $|N|\ge 30$ since the solution norm grows with the instance size. 
    \item We reformulate problem~\eqref{pr:split} as problem~\eqref{pr:reform_split}  (see Proposition~\ref{prop:slem}).
    \item To solve problem~\eqref{pr:reform_split}, we use Algorithm~\ref{alg:quadAlg} with $f_0=1\times 10^{5}$, $tol=1\times 10^{-4}$, $N=100$, $\nu_i=1$ for all $i\in [N]$.
    \item   To compute the lower bound in Algorithm~\ref{alg:quadAlg}, we use the classical SDP relaxation of the quadratic problem~\eqref{pr:reform_split},  {linearizing all non-convex quadratic terms in the control variables and obtaining an SDP constraint of the size $n_y +1$}. 
\item In the ACOPF problem, some matrices $A$ in~\eqref{eq:quadG} are negative semidefinite. We do not have to project on the corresponding equality constraints~\eqref{const:reform_split2}. Instead, we replace the equality sign by ``$\le$" and add the resulting constraints to the definition of $\A$ in~\eqref{set:A}. 
\item  {We do the posterior check in step 2 of Algorithm~\ref{alg:main} for $\hat{\zeta}=0$, i.e., we check that the obtained solution is nominally feasible.} 
\end{itemize}
Setting $\nu_i=1$ implies that we do not try to improve the objective value after finding the first feasible solution. We considered setting $\nu_i<1$, however, for all tested cases the initially obtained objective value was close to the lower bound. Attempts to improve the objective values resulted in using substantially more or all $N$ iterations with negligible or no improvement in the objective value. This result can be explained by the high quality of the lower bound from the SDP relaxation, which frequently provides a feasible or close to feasible solution to problem~\eqref{pr:reform_split}. Therefore, we decided to stop at the first obtained feasible solution by setting $\nu_i=1$ for all $i\in [N]$.

As benchmarks, we use three other approaches. First, we consider the nominal solution found by M{\sc atpower} to see how robust or not it is. Next, we implement the classical DCOPF relaxation of ACOPF, which is a linear problem. We follow the approach from~\cite{ChanceBienstock}. In this paper, the parameters $\alpha$ are optimized while in our case they are fixed. We can eliminate all equalities and second-stage variables in DCOPF, thus obtaining a standard linear robust optimization problem under ellipsoidal uncertainty. Finally, we implement an approach based on the classical SDP relaxation of ACOPF \cite{LavaeiLow}. Without the uncertainty, this relaxation usually provides exceptionally good approximations. A number of papers in the literature (e.g., \cite{venzke2017,Xie18}) implement an ARO version of this relaxation by replacing the nominal positive semidefinite matrix in the relaxation by a linear matrix function of the uncertainty; see the earlier mentioned papers for details. Since the initial relaxation is linear, except for the SDP constraint, after substituting the linear decision rule, all equalities can be eliminated and the standard linear robust approach applies. To ensure approximate feasibility of the uncertain  SDP constraint, one usually imposes this constraint for many enough realizations of the uncertainty. We followed the approach in  \cite{venzke2017} and imposed the SDP constraint for the largest possible value of each uncertainty. We have also tried using the sampling approach to ensure feasibility of the constraint with high probability, as in \cite{Xie18}. However, even for small cases this required many scenarios, resulted in too many SDP constraints, and lead to numerically unstable problems even if we used row generation to add the constraints.

For all algorithms, we limit the running time to 15 minutes (for our approach, 15 minutes from the moment after obtaining the SDP lower bound).  {DCOPF and Taylor approach do not reach this time limit in the presented cases.  For the SDP relaxation, this limit is reached for the large cases (57 and 118 buses), therefore this approach is omitted in the corresponding table. }

The two benchmark models we consider have affine constraints and use a chance constraints framework. However, to obtain the final reformulations,  both models exploit the analogy between the ellipsoidal uncertainty and chance constraints for normally distributed random variables and affine constraints. Therefore, we also use this analogy to conduct experiments.  Namely, we say that the uncertainty vector $\zeta$ is a random variable with the mean zero and the vector of standard deviations is equal to $w\%$ of the initial load. Let $\zeta$ be normally distributed, and let $W$ be a variance-covariance matrix of $\zeta$. Let $P$ be such that $W = PP\tr $ (obtained, e.g., via Cholesky decomposition) and $q_{1-\eps}$ be the $(1-\eps)$-percentile of the standard normal distribution. Then  making one affine constraint robust against the uncertainty set
\begin{align}
 \Omega:=\{\zeta: \|P^{-1} \zeta \| \le q_{1-\eps}\}   
\end{align}
 is equivalent to requiring this constraint to hold with probability at least  $ {1-}\eps\%$ for $\zeta$. This is a well-known result which follows from  classical robust optimization (see, e.g., \cite{BeTal2009}) applied to a linear constraint with ellipsoidal uncertainty and a straightforward reformulation of a chance constraint using the cumulative probability function.

We allow the uncertainty to occur at all buses with positive active power loads. This is easy to change. For instance, we could also consider uncertainty in generation by allowing additional sources of uncertainty at generator buses. We use two options for the variance-covariance matrix $W$. For the first option, we set $ W = \text{Diag }( w P^d_k)$, where $\text{Diag }$ is the operator that creates a diagonal matrix from a vector. We let $w$ vary from $0.01$ ($1\%$ of the load at the corresponding bus) to $0.5$ ($50\%$ of the load at the corresponding bus). In the second option, we allow  random correlations among the uncertainties. Namely, we generate a random positive definite matrix and scale it such that the standard deviations are equal to $w P^d_k$. We relate our uncertainty to $95\%$ chance constraints using $\Sigma= {P^{-1} }\tr P^{-1} $, $\rho =0$ and  $r=-1.65^2$ in~\eqref{def:OmegaQuad}. 
In real-life applications, one can choose a correlation matrix that is suitable for the given application; for instance, one can assume that correlations are proportional to distances between buses.

Choosing the probability of violations of the chance constraints is not trivial. The number of constraints in the problem is large even for small cases. Thus, when we restrict the probability that each chance constraint is violated, the probability that \emph{at least} one constraint is violated tends to one when the number of constraints grows. Formally, to obtain the $5\%$-safe radius of the ellipsoid for $m$ constraints, we can just require all chance constraints to hold with probability $1-\eps/m$. However, in practice this approach is usually too conservative, and thus we choose the typical $5\%$ level for each constraint and investigate how robust the solutions are for this setup.

Another reason to connect robust constraints to a probabilistic setup is the property that the volume of the ball with a constant radius tends to zero when the dimension of the ball grows. Thus, for larger cases, it would be hard to generate a large enough sample of uncertainties within the unit ball while we can easily generate uncertainties from a normal distribution.

To evaluate the performance of the algorithms, we generate $\zeta$ according to the normal distribution around zero with the covariance matrix $W$. In particular, we first generate a standard normal vector $z$ and then use $\zeta = Pz$, where $P$  comes from the earlier mentioned factorization $W = PP\tr $. We generate $1000$ values of the uncertainly including the nominal case where $\zeta=0$. 

The computational results are presented in Tables \ref{tab:table1nc}--\ref{tab:table2nc}. Tables  \ref{tab:table1nc} and \ref{tab:table1c} cover power networks with fewer than $30$ buses, and Table~\ref{tab:table2nc} shows the results for larger power networks. Our main indicator is the share of experiments where the solution provided by the corresponding approach had some infeasibilities. We round all infeasibilities down to $1\times 10^{-3}$ per unit since this is precise enough for our application. We also checked more precise rounding, and the differences in the tables were minor. In all presented instances, the equality constraints were always feasible.  

For each case, we present all levels of uncertainty where at least one model has at most $10\%$ of constraint violations.
As we mentioned earlier, it is not obvious which performance should be considered good, especially for large cases. A perfect performance would be to have less than $5\%$ of experiments with constraints violations  since this is the violation possible in our setup for a single constraint. However, given that even the smallest case we present has more than ten constraints, we choose $10\%$ as an acceptable violation.

We do not present typical M{\sc atpower}  smallest cases with three buses (LMBM3) and five buses (WB5) since for the first case even the nominal solution is robustly feasible for moderate values of the uncertainty, and for the second case none of the models is able to obtain a robustly feasible solution. The latter test case is known to be especially challenging in various contexts (see~\cite{molzahn2017}). 

In all tables, the first column denotes the values of the uncertainty as a fraction of load, mentioned earlier as $w$, in percent. The second column indicates the  evaluated model. The first two model names speak for themselves, the third name ``SDP" means the model from \cite{venzke2017}, and the name ``Taylor" indicates the approach from this paper. The third column shows average estimated objective values among all instances.  In general, robust optimization approaches do not aim to optimize the average objective, so the second column is just a rough indicator of the performance of each approach in terms of the objective value. We find this indicator sufficient since our primary goal is not to compare the objectives but rather to ensure feasibility of the constraints. The fourth column shows the running time of each benchmark algorithm required to find a solution used for the experiments. We do not include the input  construction time for each test case since all benchmarks use roughly the same inputs. The fifth column highlighted by grey shows our main indicator---the share of all experiments where the solution provided by the corresponding approach had some infeasibilities. The last two columns show the maximum number of violated constraints among all experiments: ``PQ" indicates violations of active and reactive power bounds in one test case and ``VI" indicates violations of voltage and current flow bounds in one test case. Small violations of the ``PQ" constraints could lead to more serious problems than small violations of the ``VI" constraints, therefore we have separated these indicators.

A hyphen indicates that the corresponding model delivered no solution (the resulting problem was infeasible). The bold font in the ``Model" column indicates the best model: the one which is most robust and runs in the shortest time. We have also marked by the bold font in the ``Constraint violations" column all best performing models in terms of feasibility. 

\begin{table}[!htpb]
\caption{{Results for instances up to 30 buses, \emph{without correlation}. All objective values are divided by 100 in comparison to the original data.}}
\begin{center}
     \renewcommand{\tabcolsep}{5pt}
\renewcommand{\arraystretch}{1.2} \label{tab:table1nc}
\footnotesize{
\begin{tabular}{c|l|c|c|g|c|c}   \hline 
\begin{tabular}[c]{@{}c@{}} Uncertainty, \\ \% of load  \end{tabular}
& Model   & \begin{tabular}[c]{@{}c@{}} Average \\ objective   \end{tabular}  & Time, sec. & \begin{tabular}[c]{@{}c@{}} Constraint \\  violations, \%   \end{tabular}  & \begin{tabular}[c]{@{}c@{}}Max \# violations,\\ per exper.,  PQ\end{tabular} & \begin{tabular}[c]{@{}c@{}}Max \# violations,\\ per exper., VI\end{tabular} \\ \hline
\multicolumn{ 7}{c}{case 6ww, 6 buses} \\ \hline 
\multirow{4}{*}{1}                  & Nominal & 31.3              & 0.0       & 43.9           & 0                                                                                & 2                                                                               \\
                                    & DCOPF   & -                 & -         & -              & -                                                                                & -                                                                               \\
                                    & SDP     & 31.4              & 28.4      & 22.1           & 0                                                                                & 2                                                                               \\
                                    &  \textbf{Taylor}  & 31.6              & 69.4      & \textbf{0.0}           & 0                                                                                & 0                                                                               \\ \hline
\multicolumn{ 7}{c}{case 9, 9 buses} \\ \hline 
\multirow{4}{*}{1}                  & Nominal & 53.0              & 0.0       & 0.7            & 0                                                                                & 2                                                                               \\
                                    & \textbf{DCOPF}   & 53.2              & 9.1       &  \textbf{0.0}              & 0                                                                                & 0                                                                               \\
                                    & SDP     & 53.0              & 22.5      & 100.0          & 0                                                                                & 2                                                                               \\
                                    &  Taylor  & 53.3              & 57.6      &  \textbf{0.0}              & 0                                                                                & 0                                                                               \\ \hline
\multirow{4}{*}{5}                  & Nominal & 53.2              & 0.0       & 35.7           & 0                                                                                & 2                                                                               \\
                                    & \textbf{DCOPF}   & 53.5              & 13.8      & \textbf{0.0}             & 0                                                                                & 0                                                                               \\
                                    & SDP     & 53.2              & 30.5      & 100.0          & 0                                                                                & 3                                                                               \\
                                &     Taylor  & 53.4              & 72.7      &  \textbf{0.0}              & 0                                                                                & 0                                                                               \\ \hline
\multirow{4}{*}{10}                 & Nominal & 53.5              & 0.0       & 43.9           & 0                                                                                & 2                                                                               \\
                                    & DCOPF   & 54.4              & 12.6      & 0.1            & 1                                                                                & 0                                                                               \\
                                    & SDP     & 53.6              & 23.8      & 87.5           & 0                                                                                & 3                                                                               \\
                                    & \textbf{Taylor}  & 53.7              & 70.4      & \textbf{0.0}             & 0                                                                                & 0                                                                               \\ \hline
\multirow{4}{*}{20}                 & Nominal & 54.9              & 0.0       & 48.5           & 1                                                                                & 5                                                                               \\
                                    & DCOPF   & 55.5              & 12.1      & 3.4            & 1                                                                                & 0                                                                               \\
                                    & SDP     & 55.0              & 29.4      & 99.8           & 1                                                                                & 6                                                                               \\
                                    &\textbf{Taylor}   & 55.0              & 74.4      & \textbf{1.0}            & 1                                                                                & 0                                                                               \\ \hline
\multirow{4}{*}{30}                 & Nominal & 57.1              & 0.0       & 51.3           & 1                                                                                & 6                                                                               \\
                                    & DCOPF   & -               & -      &-          &-                                                                              & -                                                                             \\
                                    & SDP     & 57.4              & 26.1      & 97.8           & 1                                                                                & 6                                                                               \\
                                     & \textbf{Taylor}   & 57.2              & 68.1      & \textbf{7.1}            & 1                                                                                & 2                                                                               \\ \hline
\multicolumn{ 7}{c}{case 30, 30 buses} \\ \hline 
\multirow{4}{*}{1}                  & Nominal & 6.1               & 0.1       & 100.0          & 0                                                                                & 2                                                                               \\
                                    & DCOPF   & -                 & -         & -              & -                                                                                & -                                                                               \\
                                    & SDP     & 5.8               & 178.6     & 31.7           & 0                                                                                & 2                                                                               \\
                                    & \textbf{Taylor}   & 5.9               & 177.2     & \textbf{1.7}            & 0                                                                                & 2                                                                               \\ \hline
\end{tabular}
}
\end{center}
\end{table}

\begin{table}[!htbp]
\caption{{Results for instances up to 30 buses, \emph{with random correlation}. All objective values are divided by 100 in comparison to the original data.}}
\begin{center}
     \renewcommand{\tabcolsep}{5pt}
\renewcommand{\arraystretch}{1.2} \label{tab:table1c}
\footnotesize{
\begin{tabular}{c|l|c|c|g|c|c}   \hline 
\begin{tabular}[c]{@{}c@{}} Uncertainty, \\ \% of load  \end{tabular}
& Model   & \begin{tabular}[c]{@{}c@{}} Average \\ objective   \end{tabular}  & Time, sec. & \begin{tabular}[c]{@{}c@{}} Constraint \\  violations, \%   \end{tabular}  & \begin{tabular}[c]{@{}c@{}}Max \# violations,\\ per exper.,  PQ\end{tabular} & \begin{tabular}[c]{@{}c@{}}Max \# violations,\\ per exper., VI\end{tabular} \\ \hline
\multicolumn{ 7}{c}{case 6ww, 6 buses} \\ \hline 
\multirow{4}{*}{1}                  & Nominal & 31.3              & 0.1       & 43.9                     & 0                                                                                & 2                                                                               \\
                                    & DCOPF   & -                 & -         & -                        & -                                                                                & -                                                                               \\
                                    & SDP     & 31.4              & 80.9      & 25.3                     & 0                                                                                & 2                                                                               \\
                                    & \textbf{Taylor}  & 31.6              & 176.9     & \textbf{0.0}                      & 0                                                                                & 0                                                                               \\ \hline
\multirow{4}{*}{5}                  & Nominal & 31.4              & 0.1       & 47.0                     & 0                                                                                & 3                                                                               \\
                                    & DCOPF   & -                 & -         & -                        & -                                                                                & -                                                                               \\
                                    & SDP    & 31.4              & 80.8      & 37.6                    & 0                                                                                & 3                                                                               \\
                                    &\textbf{Taylor}  & 32.0              & 177.0     & \textbf{4.8}                      & 0.0                                                                              & 2.0                                                                             \\ \hline
 \multicolumn{ 7}{c}{case 9, 9 buses} \\ \hline 
\multirow{4}{*}{1}                  & \textbf{Nominal} & 53.0              & 0.0       &  \textbf{0.0}                                  & 0                                                                                & 0                                                                               \\
                                    & DCOPF   & 53.2              & 33.9      &  \textbf{0.0}                                  & 0                                                                                & 0                                                                               \\
                                    & SDP     & 53.0              & 81.7      & 100.0                              & 0                                                                                & 2                                                                               \\
                                    & Taylor  & 53.3              & 201.5     &  \textbf{0.0}                                  & 0                                                                                & 0                                                                               \\ \hline
\multirow{4}{*}{5}                  & Nominal & 53.1              & 0.0       & 23.3                               & 0                                                                                & 2                                                                               \\
                                    & \textbf{DCOPF}   & 53.3              & 33.5      &  \textbf{0.0}                                  & 0                                                                                & 0                                                                               \\
                                    & SDP     & 53.1              & 81.8      & 100.0                              & 0                                                                                & 3                                                                               \\
                                    & Taylor  & 53.3              & 199.5     &  \textbf{0.0}                                  & 0                                                                                & 0                                                                               \\ \hline
\multirow{4}{*}{10}                 & Nominal & 53.2              & 0.0       & 37.6                               & 0                                                                                & 3                                                                               \\
                                    & \textbf{DCOPF}   & 53.8              & 34.2      &  \textbf{0.0}                                  & 0                                                                                & 0                                                                               \\
                                    & SDP     & 53.3              & 81.7      & 100.0                              & 0                                                                                & 4                                                                               \\
                                    & Taylor  & 53.5              & 201.6     &  \textbf{0.0}                                  & 0                                                                                & 0                                                                               \\ \hline
\multirow{4}{*}{20}                 & Nominal & 54.0              & 0.0       & 47.8                               & 1                                                                                & 5                                                                               \\
                                    & DCOPF   & 55.4              & 34.0      & 2.6                                & 1                                                                                & 0                                                                               \\
                                    & SDP     & 54.2              & 82.0      & 0.7                                & 1                                                                                & 0                                                                               \\
                                    & \textbf{Taylor}  & 54.1              & 201.0     & \textbf{0.1}                                & 1                                                                                & 0                                                                               \\ \hline
\multirow{4}{*}{30}                 & Nominal & 55.1              & 0.0       & 53.7                               & 1                                                                                & 5                                                                               \\
                                    & DCOPF   & 55.5              & 33.9      & 6.1           & 1                                                                                & 1                                                                               \\
                                    & SDP     & 55.1              & 82.3      & 68.6                               & 1                                                                                & 5                                                                               \\
                                    & \textbf{Taylor}  & 55.2              & 201.0     & \textbf{2.1}                                & 1                                                                                & 1                                                                               \\ \hline
\multirow{4}{*}{40}                 & Nominal & 56.7              & 0.0       & 57.8                               & 1                                                                                & 5                                                                               \\
                                    & DCOPF   & -                 & -         &-              & -                                                                                & -                                                                               \\
                                    & SDP     & 56.8              & 82.0      & 93.0           & 1                                                                                & 6                                                                               \\
                                    & \textbf{Taylor}  & 56.8              & 201.5     & \textbf{6.5}            & 1                                                                                & 2                                                                               \\ \hline
\multicolumn{ 7}{c}{case 30, 30 buses} \\  
\multicolumn{ 7}{c}{Infeasible for all robust approaches.} \\ \hline 
\end{tabular}
}
\end{center}
\end{table}

\begin{table}[!htbp]
\caption{Results for instances larger than 30 buses, \emph{without correlation}. All objective values are divided by 100 in comparison to the original data. \new{The SDP approach results in too-large-to-solve problems for these instances, so we do not mention this approach in the table.}}
\begin{center}
     \renewcommand{\tabcolsep}{5pt}
\renewcommand{\arraystretch}{1.2} \label{tab:table2nc}
\footnotesize{
\begin{tabular}{c|l|c|c|g|c|c}   \hline 
\begin{tabular}[c]{@{}c@{}} Uncertainty, \\ \% of load  \end{tabular}
& Model   & \begin{tabular}[c]{@{}c@{}} Average \\ objective   \end{tabular}  & Time, sec. & \begin{tabular}[c]{@{}c@{}} Constraint \\  violations, \%   \end{tabular}  & \begin{tabular}[c]{@{}c@{}}Max \# violations,\\ per exper.,  PQ\end{tabular} & \begin{tabular}[c]{@{}c@{}}Max \# violations,\\ per exper., VI\end{tabular} \\ \hline
\multicolumn{ 7}{c}{case 57, 57 buses} \\ \hline 
\multirow{4}{*}{1}                  & Nominal & 417.4             & 0.1       & 70.6           & 2                                                                                & 1                                                                               \\
                                    & DCOPF   & 418.5             & 43.0      & 100.0          & 2                                                                                & 1                                                                               \\
                                    & \textbf{Taylor}  & 426.8             & 467.2     & \textbf{0.0}            & 0                                                                                & 0                                                                               \\ \hline
\multicolumn{ 7}{c}{case 118, 118 buses} \\ \hline 
\multirow{4}{*}{1}                  & Nominal & 1296.7            & 0.2       & 99.5           & 9                                                                                & 0                                                                               \\
                                    & DCOPF   & 1315.6            & 94.9      & 100.0          & 21                                                                               & 0                                                                               \\
                                    &  \textbf{Taylor}   & 1301.3            & 830.0     &\textbf{1.1}            & 1                                                                                & 0                                                                               \\ \hline
\end{tabular} \\ \medskip 
}
\end{center}
\end{table}

In most cases,  the nominal solution is not robustly feasible, thus considering the uncertainty is important. The results show that our approach provides a robust solution in reasonable time even for the instance with 118 buses. 
 By Corollary~\ref{cor:main}, the approach we implemented is guaranteed to be robust only for small uncertainty sets, which is indeed reflected in the numerical results. For all presented cases the ``Taylor" solution is robust in more than $95\%$ of the experiments when the uncertainty is $1\%$ of the loads while this is not the case for the nominal solutions and other approaches. The second-best approach in terms of robustness is ``DCOPF". It works for some small instances but fails for larger ones. DCOPF is an accurate approximation under quite restrictive assumptions and seems to become less precise under uncertainty. Finally the ``SDP" model provides most violations among the three robust models.  However, if we look deeper into these violations, when the uncertainty sets are small, the ``SDP"  solution usually leads to few small violations of ``VI" constraints. This could be expected since the voltages are  approximated in the SDP relaxation and these approximations could become looser under uncertainty. Adding some tightening constraints to the ``SDP" model could potentially improve the robustness of that approach. 

As to the running times, our approach is the slowest for small test cases since it is an iterative approach due to Algorithm~\ref{alg:quadAlg}. However, when our approach found robust  solutions, it could usually do that within several iterations of the alternating projections Algorithm~\ref{alg:quadAlg}. Longer running times for the two largest instances indicate that the alternating projections algorithm went through many iterations and could not converge, thus a robust solution may not exist for the considered subset of state variables. For the two largest cases, the SDP approach is the slowest, and it is in fact too large to solve (e.g., for case 118 it includes $100$ SDP constraints of the size $236\times 236$).

Finally, our approach often provides somewhat higher average objective than others.  Given that the approach is more robust than the others, the difference in the average objective is minimal. To understand the performance with respect to the objective better, we looked at Case 9 with correlations, which is robustly feasible for three of four models for the $1\%$-uncertainty and robustly feasible for two models for the $10\%$-uncertainty. Below we show the box plots of the objective values of all robustly feasible approaches. The two box plots use same interval of possible objective values in the $y$-axes for comparability.  
\begin{figure}[!htbp]%
\caption{ {Box plot of the objective value realizations for models without constraint violations for Case 9 with correlations,  $1\%$ and $10\%$  uncertainty.  \medskip }}
\centering
\includegraphics[width=0.5\textwidth]{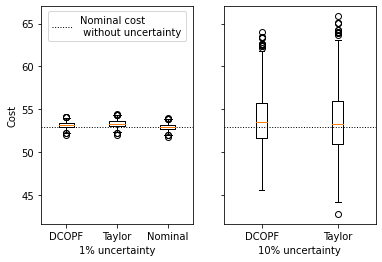}
\end{figure} 
We see that the ``Taylor" approach tends to have higher cost in both cases, and it also results in larger standard deviations for the larger uncertainty, meaning that very low costs are also more common than in DCOPF. The worst-case objective value of our algorithm in the experiments is a little larger than for DCOPF, and the algorithm does not find the DCOPF solution. To find more solutions, such as the DCOPF solution when it is robust, it could be useful to enlarge the search space of the algorithm.

Finally, we see that adding correlations to the uncertainty ellipsoid does not substantially change the results for small cases. However, when the instance size grows, no robust solutions could be found by any model for ``case30", ``case57", and ``case118", even for the lowest values of the uncertainty.  Therefore, we do not provide results with correlations for larger cases. Further investigation is needed to see why exactly this situation occurs and to possibly tune the models to work for larger cases with correlated uncertainties.

Our evaluations are equally conservative for all models. First, we say that more than $10\%$ of experiments with constraint violations is a bad performance while we have many constraints  and the analogy between the normal distribution and our ellipsoid is only valid for one constraint. The $5\%$-safe radius of the ellipsoid for many constraints would be much larger, but using this radius would be too conservative. Second, we do not deeply analyze magnitudes of violations. By the experiment setup they are larger than $1\times 10^{-3}$ per unit, but the importance of such violations can depend on the type of constraint and the details of the test case. Finally, we did not fine-tune the algorithm to each specific case and allowed for as many power injection uncertainties as possible. The performance could improve if the algorithms are fine-tuned, especially for large instances.

\section{Conclusions and directions for future research} \label{sec:aviol}

In this paper, we propose a framework to obtain approximately feasible solutions to quadratic ARO problems with equalities, which implicitly define the second-stage decision rules for the state variables. We replace the implicit decision rules by their explicit piecewise affine approximations. As a result, we can eliminate the state variables from the problem and replace the original ARO problem by a sequence of classical quadratic problems with additional tractable conic constraints. Since, a generally efficient algorithm does not exist yet for the latter problems, we design an alternating projections algorithm that converges to a local optimum of the problem. For any $\eps>0$, if the piecewise affine approximations are fine enough, the suggested algorithm  guarantees that  {the second-stage equality constraints are satisfied and the inequality constraints are not violated by more than $\eps$ on ``large" subsets of uncertainty, where ``large" is defined in Theorem~\ref{thm:Taylor}~\ref{thm:Taylor1},~\ref{thm:Taylor_add}}. The feasibility of the second-stage equalities is rarely addressed in the literature, thus we consider analysing and ensuring such feasibility  an important contribution to the existing research. We suggest two versions of the algorithm, Algorithm~\ref{alg:Taylor} for general problems and Algorithm~\ref{alg:main} for problems linear in the control variables.  

We implement the algorithm for ACOPF problems with uncertainty in loads and simulate the uncertainty to evaluate the performance of the algorithm in comparison to three known benchmarks: nominal solution, robust DCOPF and robust SDP relaxation. The solutions provided by our approach are robustly feasible  for small uncertainty sets and for cases with up to 118 buses. Moreover, the solutions are substantially more robust than the benchmarks. The algorithm also performs well in terms of the objective function values. The experiments show that DCOPF and SDP approximations work well for problems without uncertainty, but possible inexactness of these methods  is magnified after adding the uncertainty. As a result, they are less robust in an ex-post experiment. In contrast, our approach yields more robust results in the experiments. The good performance can be explained by the fact that our approach preserves more non-linearities from the original problem. Notice that the size of the final problem (i.e., problem~\eqref{pr:split}) which we solve depends on the number of sources of uncertainty. Hence, even for larger networks, the approach could work well for systems with uncertain power injections at a limited subset of the buses.

We conclude by analyzing the limitations of our algorithms and suggesting directions for further research. First, the non-linearities which we keep in the problem have drawbacks. Namely, when DCOPF can find a robust solution, it does so substantially faster than the approach we propose. At the last step of our approach, after removing the state and uncertainty variables, we need to solve a problem with quadratic and SDP constraints. We proposed the alternating projections in Algorithm~\ref{alg:quadAlg} for this purpose. This algorithm only  finds locally optimal solutions, may take a number of iterations to converge, and each iteration solves an SDP problem, which makes the procedure rather slow for larger instances. However, this approach is still relevant since no generally efficient algorithms exist for such problems yet, and our Algorithm~\ref{alg:quadAlg}  contributes to the development of quadratic optimization with SDP constraints. In the future, one could use an alternative algorithm or a relaxation to solve the obtained problem with quadratic and SDP constraints.

Next, our algorithms can exploit various subsets of uncertainty, control and/or state variables. However, the actual implementation considered small uncertainty sets without partitioning them, and we restricted the search to one subset of state variables around the initial optimal solution. As expected, for larger values of the uncertainty in our numerical experiments, the algorithm becomes imprecise or could not find a feasible solution. 
A natural remedy to increase precision and find additional feasible solutions would be to consider several subsets of the state variables as described in subproblem~\eqref{pr:split}. Next, to work with larger uncertainty sets, we could combine Algorithm~\ref{alg:Taylor} 
 and Algorithm~\ref{alg:main} with the approach proposed in~\cite{splittingUnc} to split larger uncertainty sets. 

Third, some improvements could be done in particular for ACOPF problems. These problems possess much sparsity, and it is to a large extent preserved under the transformation~\eqref{eq:TaylImpl}.  Therefore, it would be also important to explore the true scalability of our method for ACOPF problems by implementing the algorithms more efficiently and investigating possibilities to exploit the inherent sparsity structure of ACOPF problems. Finally, it would be interesting to look at the extensions of our method for multiperiod problems, especially  those incorporating binary variables to enable application to unit commitment problems for longer-term generator scheduling. 

 \section*{Acknowledgment}
Olga Kuryatnikova and Bissan Ghaddar were supported by NSERC Discovery Grant 2017-04185. Daniel Molzahn was supported by the National Science Foundation under grant number 2023140. We also greatly thank the referees, the Associate Editor, and the Editor for their thorough and thoughtful comments that helped us improve the quality of the paper.

\pagebreak
\bibliographystyle{siam}

\bibliography{bibOPF}


\end{document}